\documentclass[12pt]{article}
\usepackage{amsmath,amssymb}

\newcommand{\ncm}{\newcommand}
\ncm{\aut}{auto\-mor\-phi\-sm} \ncm{\Inn}{\mbox{\rm Inn}} 
\ncm{\Ap}{\mbox{$\overline{\rm Inn}$}} \ncm{\Ext}{\mbox{\rm Ext}} 
\ncm{\Ex}{\mbox{\rm Ex}} \ncm{\OExt}{\mbox{\rm OrderExt}} 
\ncm{\AI}{\mbox{\rm AInn}} \ncm{\HI}{\mbox{\rm HInn($A$)}} 
\ncm{\Aut}{\mbox{\rm Aut}} \ncm{\Mal}{\mbox{$M_{\alpha}$}} 
\ncm{\Aff}{\mbox{${\rm Aff}$}} \ncm{\id}{\mbox{\rm id}} 
\ncm{\Ker}{\mbox{\rm Ker}} \ncm{\BE}{\begin{eqnarray*}} 
\ncm{\EE}{\end{eqnarray*}} \ncm{\lra}{\mbox{$\longrightarrow$}} 
\ncm{\Hom}{\mbox{\rm Hom}} \ncm{\calU}{{\cal U}} \ncm{\el}{\ell} 
\ncm{\ad}{\mbox{\rm ad}} \ncm{\Alg}{\mbox{\rm Alg}} 
\ncm{\Conv}{\mbox{\rm Conv}} \ncm{\D}{{\cal D}} 

\ncm{\cstar}{$C^{*}$-algebra} \ncm{\cstars}{$C^{*}$-algebras} 
\ncm{\ra}{\mbox{$\rightarrow$}} \ncm{\la}{\mbox{$\leftarrow$}} 
\ncm{\hra}{\hookrightarrow} \ncm{\da}{\mbox{$\downarrow$}} 
\ncm{\se}{\mbox{$\searrow$}} \ncm{\al}{\mbox{$\alpha $}} 
\ncm{\del}{\mbox{$\delta$}} \ncm{\supp}{\mbox{\rm supp}} 
\ncm{\Ad}{\mbox{\rm Ad}} \ncm{\CAR}{\mbox{$M_{2^{\infty}}$}} 
\ncm{\ep}{\mbox{$\epsilon > 0$}} \ncm{\ol}{\overline} 
\ncm{\Mninf}{\mbox{$M_{n^{\infty}}$}} \ncm{\MR}{M. R\o{}rdam} 
\ncm{\Range}{\mbox{\rm Range}} 
\ncm{\vo}{}
\ncm{\ch}{}
\ncm{\CMP}{Comm. Math. Phys.} \ncm{\add}{} 
\ncm{\tilsig}{\tilde{\sigma}} \ncm{\dist}{{\rm 
dist}}\ncm{\eps}{\epsilon}  \ncm{\calL}{{\mathcal{L}}}
 
\ncm{\calH}{{\mathcal{H}}} 
\ncm{\lan}{{\langle}}\ncm{\ran}{{\rangle}}

\newtheorem{theo}{Theorem}[section]

\newtheorem{lem}[theo]{Lemma}
\newtheorem{prop}[theo]{Proposition}
\newtheorem{remark}[theo]{Remark}
\newtheorem{definition}[theo]{Definition}
\newtheorem{example}[theo]{Example}
\newtheorem{property}[theo]{Property}

\newenvironment{rem}{\begin{remark} \rm}{\end{remark}}
\newenvironment{pf}{{\it Proof.}}{\hfill$\square$\vspace{3mm}}

\ncm{\R}{\mbox{\bf R}} \ncm{\Z}{\mbox{\bf Z}} \ncm{\T}{\mbox{\bf 
T}} \ncm{\TT}{\T$^{2}$} \ncm{\N}{\mbox{\bf N}} \ncm{\C}{\mbox{\bf 
C}} 

\title{Homogeneity of the pure state space for separable \cstars}
\bigskip
\author{Hajime Futamura, 
Nobuhiro Kataoka, and Akitaka Kishimoto\\
 Department of Mathematics, Hokkaido University,
 Sapporo 060, Japan}
\date{January 2001}
\begin{document}
\maketitle 

\begin{abstract} We prove that the pure state space is 
homogeneous under the action of the automorphism group (or a 
certain smaller group of approximately inner automorphisms) for a 
fairly large class of simple separable nuclear \cstars, including 
the approximately homogeneous \cstars\ and the class of purely 
infinite \cstars\ which has been recently classified by Kirchberg 
and Phillips. This extends the known results for UHF algebras and 
AF algebras by Powers and Bratteli. \\

 \hfill Mathematics Subject 
Classification: 46L40, 46L30 
\end{abstract}

 \ncm{\U}{{\mathcal{U}}} 
 \ncm{\F}{{\mathcal{F}}}
 \ncm{\G}{{\mathcal{G}}} 
 \ncm{\om}{{\omega}}
 \ncm{\Hil}{{\mathcal{H}}}
    
\section{Introduction}
If $A$ is a \cstar, we denote by $S(A)$ the convex set of states 
of $A$ and by $P(A)$ the set of pure states of $A$. The 
automorphism group $\Aut(A)$ of $A$ induces an action $\phi$ on 
$S(A)$ by affine homeomorphisms, i.e., if $\alpha\in\Aut(A)$, then 
$\phi(\alpha)$ sends $f\in S(A)$ to $f\alpha^{-1}$. Note that this 
action leaves $P(A)$ invariant since $P(A)$ is the set of extreme 
points of $S(A)$ and hence that it cannot be transitive on $S(A)$ 
(except for the trivial case $A=\C1$). But from Powers' result 
\cite{Pow} for UHF algebras we know that $\Aut(A)$ can act 
transitively on $P(A)$ for some simple \cstars. See \cite{Br} for 
an immediate extension to AF algebras and \cite{BK1} for a partial 
extension to Cuntz algebras. 

Note that $\Aut(A)$ has some distinguished subgroups; $\Inn(A)$, 
$\AI(A)$, $\Ap(A)$. Denoting by $\U(A)$ the unitary group of $A$ 
(or $A+\C1$ if $A\not\ni1$), the group $\Inn(A)$ of inner 
automorphisms is given by $\{\Ad\,u\ |\ u\in\U(A)\}$ and the group 
$\Ap(A)$ of approximately inner automorphisms is the closure of 
$\Inn(A)$ in $\Aut(A)$ with the topology of strong convergence. 
The group $\AI(A)$ of asymptotically inner automorphisms consists 
of $\alpha\in\Ap(A)$ such that there exists a continuous path 
$(u_t)_{t\in [0,1)}$ in $\U(A)$ with 
$\alpha=\lim_{t\rightarrow1}\Ad\,u_t$; in the classification 
theory of purely infinite simple \cstars\ \cite{KP}, $\AI(A)$ is 
characterized as the group of automorphisms which have the same KK 
class as the identity automorphism (see \cite{KK} for a similar 
characterization in the case of simple unital AT algebras of real 
rank zero). Note that $\Inn(A)\subset \AI(A)\subset \Ap(A)\subset 
\Aut(A)$ and that the inclusions are proper in general. In the 
statements given in the previous paragraph $\Aut(A)$ should be 
replaced by $\Ap(A)$ from the way in which they are proven. 

In this paper we aim to prove the following statement for a large 
class of separable \cstars\ $A$: If $\om_1$ and $\om_2$ are pure 
states of $A$ with $\ker\pi_{\om_1}=\ker\pi_{\om_2}$, then there 
is an $\alpha\in\Ap(A)$ such that $\om_1=\om_2\alpha$. Looking at 
the proofs closely, even more is true in the cases we will handle; 
it will follow that one can choose $\alpha$ from $\AI(A)$. Hence, 
restricted to the case that the \cstar\ $A$ is simple, we are to 
try to prove that $\AI(A)$ acts on $P(A)$ transitively. 

In the subsequent section we will introduce some properties for 
\cstars\ and show that these properties imply the homogeneity of 
the pure state space in the above sense if the \cstars\ are 
separable.  This section contains the main idea of this paper. 

In Sections 3--5 we will prove the above properties (more 
precisely, Property \ref{A10} below, the strongest among them) for 
a large class of \cstars, namely approximately homogeneous 
\cstars, simple crossed products of AF algebras by $\Z$, and a 
class of purely infinite \cstars\ including the class which is 
classified by Kirchberg and Phillips \cite{KP,Phi}. Hence the 
above-mentioned transitivity follows for these \cstars. We might 
expect that this transitivity holds for all separable nuclear 
\cstars. In section 6 we show that the transitivity holds also for 
the group \cstars\ if the group is a countable discrete amenable 
group.  

In Section 7 we note that even a stronger form of transitivity 
holds for the above mentioned separable \cstars. One consequence 
is that given any sequence $(\pi_n)$, indexed by $\Z$, of distinct 
points in the set of equivalence classes of irreducible 
representations there is an asymptotically inner automorphism 
$\alpha$ which shifts this sequence, i.e.,  
$\pi_n\alpha=\pi_{n+1}$ for all $n$. 

The properties for \cstars\ introduced in Section 2 seem 
interesting for their own sake. In Section 8 we will give yet 
another version, which is valid for the \cstars\ treated above, to 
derive a closure property saying that if the property holds for a 
\cstar, then it also holds for all its hereditary 
$C^*$-subalgebras. See \ref{R4} and \ref{R5} for other 
consequences. 

The third named author would like to thank O. Bratteli for 
discussions at an early stage of this work.  

\section{Homogeneity of the pure state space}
Let $A$ be a \cstar\ and $\U(A)$ the unitary group of $A$ (or 
$A+\C1$ if $A \not\ni 1$). If two pure states $\omega_1$ and 
$\omega_2$ are equivalent or $\omega_1\sim\omega_2$, i.e., the GNS 
representation $\pi_{\omega_1}$ and $\pi_{\omega_2}$ are 
equivalent, then there is a $u\in\U(A)$ such that 
$\omega_1=\omega_2\Ad\,u$ by Kadison's transitivity (1.21.16 of 
\cite{Sak}); we shall repeatedly use this fact below. First we 
define the following two properties for $A$: 

\begin{property}\label{A1} 
For any finite subset $\F$ of $A$ and $\eps>0$ there exist a 
finite subset $\G$ of $A$ and $\delta>0$ satisfying: If $\omega_1$ 
and $\omega_2$ are pure states of $A$ such that $\omega_1\sim 
\omega_2$ and 
 $$ |\omega_1(x)-\omega_2(x)|<\delta,\ \ x\in\G,
 $$
then there is a $u\in\U(A)$ such that $\omega_1=\omega_2\Ad\,u$ 
and 
 $$ \|\Ad\,u(x)-x\|<\eps,\ \ x\in\F.
 $$
\end{property} 
  
\begin{property}\label{A2} 
For any finite subset $\F$ of $A$ and $\eps>0$ there exist a 
finite subset $\G$ of $A$ and $\delta>0$ satisfying: If $\omega_1$ 
and $\omega_2$ are pure states of $A$ such that $\ker 
\pi_{\omega_1}=\ker\pi_{\omega_2}$ and 
 $$ |\omega_1(x)-\omega_2(x)|<\delta,\ \ x\in\G,
 $$
then for any finite subset $\F'$ of $A$ and $\eps'>0$ there is a 
$u\in\U(A)$ such that  
 \BE |\omega_1(x)-\omega_2\Ad\,u(x)|&<&\eps',\ \ x\in\F',\\
 \|\Ad\,u(x)-x\|&<&\eps,\ \ x\in\F.
 \EE
\end{property} 

The following lemma is well-known (cf. \cite{BKR}).

\begin{lem}\label{A3}
Let $\om_1$ and $\om_2$ be pure states of $A$ such that 
$\ker\pi_{\om_1}=\ker\pi_{\om_2}$. For any finite subset $\F$ of 
$A$ and $\eps>0$ there is a $u\in\U(A)$ such that 
 $$
 |\om_1(x)-\om_2\Ad\,u(x)|<\eps,\ \ x\in\F.
 $$
\end{lem}
\begin{pf}
Since $\om_1$ is a pure state, there exists an $e\in A$ such that 
$e\geq0,\, \|e\|=1,\, \om_1(e)=1$, and 
 $$
 \|exe-\om_1(x)e^2\|<\eps,\ \ x\in\F.
 $$
Since $\|\pi_{\om_1}(e)\|=1$, we may suppose, by slightly changing 
$e$ if necessary, that there is a $\xi\in \Hil_{\pi_{\om_2}}$ such 
that $\|\xi\|=1$ and $\pi_{\om_2}(e)\xi=\xi$. Since $\pi_{\om_2}$ 
is irreducible there is a $u\in\U(A)$ such that 
$\pi_{\om_2}(u^*)\Omega_{\om_2}=\xi$. Since $\langle 
\pi_{\om_2}(exe)\xi,\xi\rangle=\om_2\Ad\,u(x)$, we obtain the 
conclusion. \end{pf}                  

\begin{prop} \label{A4}
For any \cstar\ $A$ Property \ref{A1} implies Property \ref{A2}.
\end{prop}
\begin{pf}
Suppose \ref{A1} and choose $(\G,\delta)$ for $(\F,\eps)$ as 
therein. Suppose that $\om_1$ and $\om_2$ are pure states of $A$ 
such that $\ker\pi_{\om_1}=\ker\pi_{\om_2}$ and
 $$ |\om_1(x)-\om_2(x)|<\delta/2,\ \ x\in\G.
 $$
Let $\F'$ be a finite subset of $A$ and $\eps'>0$ with 
$\eps'<\delta/2$. Then by Lemma \ref{A3} there is a $v\in\U(A)$ 
such that 
 $$
 |\om_2\Ad\,v(x)-\om_1(x)|<\eps',\ \ x\in\F'\cup\G.
 $$
Since $|\om_2\Ad\,v(x)-\om_2(x)|<\delta,\ x\in\G$, there is a 
$u\in\U(A)$ such that $\|\Ad\,u(x)-x\|<\eps,\ x\in\F$ and 
$\om_2\Ad\,v=\om_2\Ad\,u$. The latter condition implies that 
$|\om_2\Ad\,u(x)-\om_1(x)|<\eps',\ x\in\F'$. This completes the 
proof. 
\end{pf}  

\begin{theo}\label{A5}
Let $A$ be a separable \cstar. Then the following conditions are 
equivalent: 
 \begin{enumerate}
 \item  Property \ref{A2} holds.
 \item  For any finite subset $\F$ of $A$ and $\eps>0$ there exist a 
finite subset $\G$ of $A$ and $\delta>0$ satisfying: If $\om_1$ 
and $\om_2$ are pure states of $A$ such that 
$\ker\pi_{\om_1}=\ker\pi_{\om_2}$ and 
$|\om_1(x)-\om_2(x)|<\delta,\ x\in\G$, there is an 
$\alpha\in\Ap(A)$ such that $\om_1=\om_2\alpha$ and 
$\|\alpha(x)-x\|<\eps,\ x\in\F$.  If $\F=\emptyset$, then $\G$ can 
be assumed to be empty. 
 \end{enumerate} 
\end{theo} 
\begin{pf}          
First suppose (2). For any $(\F,\eps)$ we choose $(\G,\delta)$ as 
in (2). Let $\om_1$ and $\om_2$ be pure states of $A$ such that 
$\ker\pi_{\om_1}=\ker\pi_{\om_2}$ and 
$|\om_1(x)-\om_2(x)|<\delta,\ x\in\G$. Then there is an 
$\alpha\in\Ap(A)$ as in (2). Since there is a sequence $(u_n)$ in 
$\U(A)$ such that $\alpha=\lim\Ad\,u_n$, it follows that for any 
finite subset $\F'$ of $A$ and $\eps'>0$ there is an $n$ such that 
$|\om_1(x)-\om_2\Ad\,u_n(x)|<\eps',\ x\in\F'$ and 
$\|\Ad\,u_n(x)-x\|<\eps,\ x\in\F$. Thus (1) follows. 

Suppose (1).  Given $(\F,\eps)$ choose $(\G,\delta)$ as in 
Property \ref{A2} and let $\om_1$ and $\om_2$ be pure states of 
$A$ such that $\ker\pi_{\om_1}=\ker\pi_{\om_2}$ and  
$|\om_1(x)-\om_2(x)|<\delta,\ x\in\G$. (If $\F=\emptyset$, then we 
set $\G=\emptyset$, which will take care of the last statement.)

Let $(x_n)$ be a dense sequence in  $A$.

Let $\F_1=\F\cup\{x_1\}$ and let $(\G_1,\delta_1)$ be the 
$(\G,\delta)$ for $(\F_1,\eps/2)$ as in Property \ref{A2}.  By  
\ref{A2} there is a $u_1\in\U(A)$ such that 
 \BE 
 \|\Ad\,u_1(x)-x\|&<&\eps,\ \ x\in\F,\\ 
 |\om_1(x)-\om_2\Ad\,u_1(x)|&<&\delta_1,\ \ x\in\G_1.
 \EE
Let $\F_2=\F\cup\{x_1,x_2,\Ad\,u_1^*(x_1),\Ad\,u_1^*(x_2)\}$ and 
let $(\G_2,\delta_2)$ be the $(\G,\delta)$ for $(\F_2,\eps/2^2)$ 
as in Property \ref{A2} such that $\G_2\supset \G_1$ and 
$\delta_2<\delta_1$. Then there is a $u_2\in\U(A)$ such that 
 \BE
 \|\Ad\,u_2(x)-x\|&<&2^{-1}\eps,\ \ x\in\F_1, \\
 |\om_2\Ad\,u_1(x)-\om_1\Ad\,u_2(x)|&<&\delta_2, \ \ x\in\G_2.
 \EE
Let $\F_3=\F\cup\{x_i,\Ad\,u_2^*(x_i)\ |\ i=1,2,3\}$ and let 
$(\G_3,\delta_3)$ be the $(\G,\delta)$ for $(\F_3,2^{-3}\eps)$ as 
in Property \ref{A2} such that $\G_3\supset \G_2$ and 
$\delta_3<\delta_2$. Then there is a $u_3\in\U(A)$ such that
 \BE
 \|\Ad\,u_3(x)-x\|&<&2^{-2}\eps,\ \ x\in\F_2,\\
 |\om_1\Ad\,u_2(x)-\om_2\Ad\,u_1u_3(x)|&<&\delta_3,\ \ x\in\G_3.
 \EE
We shall repeat this process.

If $\F_k,\G_k,\delta_k,$ and $u_k$ are given for $k<n$, let
 $$
 \F_n=\F\cup\{x_i,\Ad\,u_{n-1}^*u_{n-3}^*\cdots u_{\#}^*(x_i)\ |\ 
 i=1,2,\ldots,n\},
 $$        
where $\#=2$ or $1$ depending on the parity of $n$. We then choose 
$(\G_n,\delta_n)$ as in Property \ref{A2} for $(\F_n,2^{-n}\eps)$
such that $\G_n\supset \G_{n-1}$ and $\delta_n<\delta_{n-1}$. 
Finally we pick up a $u_n\in\U(A)$ such that
 $$
 \|\Ad\,u_n(x)-x\|<2^{-n+1}\eps,\ \ x\in\F_{n-1}
 $$
and if $n$ is odd,
 $$
 |\om_1\Ad(u_2u_4\cdots u_{n-1})(x)-\om_2\Ad(u_1u_3\cdots u_n)(x)|<\delta_n, \ \ x\in\G_n,
 $$
else if $n$ is even,
 $$   
 |\om_1\Ad(u_2u_4\cdots u_{n})(x)-\om_2\Ad(u_1u_3\cdots u_{n-1})(x)|<\delta_n, \ \ 
 x\in\G_n.
 $$       
We may assume that $\cup_n\G_n$ is dense in $A$. 

Since $\|\Ad\,u_k(x_i)-x_i\|<2^{-k+1}\eps$ for $k>i$, we have that
 $$
 \lim_{n\ra\infty}\Ad(u_1u_3\cdots u_{2n-1})(x_i)
 $$
converges for any $i$. Since $(x_i)$ is a dense sequence in $A$, 
we have that $\Ad(u_1u_3\cdots u_{2n-1})$ converges strongly on 
$A$; thus the limit $\alpha$ exists as an endomorphism of $A$. In 
a similar way $\Ad(u_2u_4\cdots u_{2n})$ converges to an 
endomorphism $\beta$ of $A$. From the estimate on 
$\om_1\Ad(u_2u_4\cdots u_{2n})-\om_2\Ad(u_1u_3\cdots u_{2n+1})$, 
it follows that $\om_1\beta=\om_2\alpha$. 

We shall show that $\alpha$ and $\beta$ are automorphisms. For 
that purpose it suffices to show that 
$\Ad(u_{2n}^*u_{2n-2}^*\cdots u_2^*)$ and a similar expression 
with odd-numbered $u_k$ converge as $n\ra\infty$. 

Since $\|\Ad\,u_n(x)-x\|<2^{-n+1}\eps,\ x\in\F_{n-1}$, we have 
that if $n>i$,
 $$
 \|\Ad\,u_n\,\Ad(u_{n-2}^*u_{n-4}^*\cdots u_{\#}^*)(x_i)-  
 \Ad(u_{n-2}^*u_{n-4}^*\cdots u_{\#}^*)(x_i)\|<2^{-n+1}\eps.
 $$
This implies the desired convergence. 

Finally, since $\F_n\supset\F$, it follows that the automorphisms 
$\alpha=\lim\Ad(u_1u_3\cdots u_{2n-1})$ and 
$\beta=\lim\Ad(u_2u_4\cdots u_{2n})$ satisfy that 
 $$\|\alpha(x)-x\|<4\eps/3,\ \ \|\beta(x)-x\|<2\eps/3,
 $$
for $x\in\F$. Hence $\|\alpha\beta^{-1}(x)-x\|<2\eps,\ x\in\F$. 
Since $\om_1=\om_2\alpha\beta^{-1}$, this completes the proof. 
\end{pf}
   
In the following sections we will actually treat properties 
stronger than \ref{A1}. 

\begin{property}\label{A7}
For any finite subset $\F$ of $A$ and $\eps>0$ there exist a 
finite subset $\G$ of $A$ and $\delta>0$ satisfying: If $\omega_1$ 
and $\omega_2$ are pure states of $A$ such that $\omega_1\sim 
\omega_2$ and 
 $$ |\omega_1(x)-\omega_2(x)|<\delta,\ \ x\in\G,
 $$
then there is a continuous path $(u_t)_{t\in[0,1]}$ in $\U(A)$ 
such that $u_0=1$, $\omega_1=\omega_2\Ad\,u_1$, and 
 $$ \|\Ad\,u_t(x)-x\|<\eps,\ \ x\in\F,\ t\in [0,1].
 $$
\end{property}

As in the proof of \ref{A4} we can show that the above property 
implies:
  
\begin{property}\label{A8} 
For any finite subset $\F$ of $A$ and $\eps>0$ there exist a 
finite subset $\G$ of $A$ and $\delta>0$ satisfying: If $\omega_1$ 
and $\omega_2$ are pure states of $A$ such that $\ker 
\pi_{\omega_1}=\ker\pi_{\omega_2}$ and 
 $$ |\omega_1(x)-\omega_2(x)|<\delta,\ \ x\in\G,
 $$
then for any finite subset $\F'$ of $A$ and $\eps'>0$ there is a 
continuous path $(u_t)_{t\in[0,1]}$ in $\U(A)$ such that $u_0=1$, 
and 
 \BE |\omega_1(x)-\omega_2\Ad\,u_1(x)|&<&\eps',\ \ x\in\F',\\
 \|\Ad\,u_t(x)-x\|&<&\eps,\ \ x\in\F,\ t\in[0,1].
 \EE
\end{property} 

We recall here that $\AI(A)$ is the group of asymptotically inner 
automorphisms of $A$, a proper normal subgroup of $\Ap(A)$ in 
general. The following result can be shown exactly in the same way 
as Theorem \ref{A5} is shown:

\begin{theo} \label{A9}
Let $A$ be a separable \cstar. Then the following conditions are 
equivalent:
 \begin{enumerate}
 \item Property \ref{A8} holds.
 \item For any finite subset $\F$ of $A$ and $\eps>0$ there exist a 
finite subset $\G$ of $A$ and $\delta>0$ satisfying: If $\om_1$ 
and $\om_2$ are pure states of $A$ such that 
$\ker\pi_{\om_1}=\ker\pi_{\om_2}$ and 
$|\om_1(x)-\om_2(x)|<\delta,\ x\in\G$, there exist an 
$\alpha\in\AI(A)$ and a continuous path $(u_t)_{t\in[0,1)}$ in 
$\U(A)$ such that $u_0=1$, $\alpha=\lim_{t\rightarrow1}\Ad\,u_t$, 
$\om_1=\om_2\alpha$ and $\|\Ad\,u_t(x)-x\|<\eps,\ x\in\F,\, t\in 
[0,1)$.  If $\F=\emptyset$, then $\G$ can be assumed to be empty.  
 \end{enumerate}
\end{theo}
\begin{pf} 
The proof proceeds exactly in the same way as the proof of  
\ref{A5}. In the proof of (1)$\Rightarrow$(2) of \ref{A5} we have 
defined the automorphisms $\alpha$ and $\beta$ of $A$; the 
$\alpha$ in the above statement is $\alpha\beta^{-1}$. The 
$\alpha$ in the proof of \ref{A5} is defined as the limit of 
$\Ad(u_1u_3\cdots u_{2n-1})$. In the present assumption we have a 
continuous path $(u_{nt})$ in $\U(A)$ for each $n$ such that 
$u_{n0}=1, \ u_{n1}=u_n$, and $\|\Ad\,u_{nt}(x)-x\|<2^{-n+1}\eps,\ 
x\in\F_{n-1}$. We define a continuous path 
$(v_t)_{t\in[0,\infty)}$ by: for $t\in [n,n+1]$, 
  $$
  v_t=u_1 u_3\cdots u_{2n-1}u_{2n+1,t-n}.
  $$ 
Then it follows that $v_0=1$,
$\alpha=\lim_{t\rightarrow\infty}\Ad\,v_t$, and 
$\|\Ad\,v_t(x)-x\|\leq 4\eps/3,\ x\in\F$. We prove that $\beta$ 
also enjoys a similar property and thus $\alpha\beta^{-1}$ does 
too.  This shows that (1)$\Rightarrow$(2). 

The other implication is obvious since it is assumed that $u_0=1$ 
for the path $(u_t)$ in (2).
\end{pf}

We will consider even a stronger property:

\begin{property}\label{A10}
For any finite subset $\F$ of $A$ and $\eps>0$ there is a finite 
subset $\G$ of $A$ and $\delta>0$ satisfying: If $B$ is a \cstar\ 
containing $A$ as a $C^*$-subalgebra and  $\omega_1$ and 
$\omega_2$ are pure states of $B$ such that $\omega_1\sim 
\omega_2$ and 
 $$ |\omega_1(x)-\omega_2(x)|<\delta,\ \ x\in\G,
 $$
then there is a continuous path $(u_t)_{t\in[0,1]}$ in $\U(B)$ 
such that $u_0=1$, $\omega_1=\omega_2\Ad\,u_1$, and 
 $$ \|\Ad\,u_t(x)-x\|<\eps,\ \ x\in\F,\ t\in [0,1].
 $$
\end{property}

In the above property we may replace the exact equality 
$\omega_1=\omega_2\Ad\,u_1$ by an approximate one
$\|\omega_1-\omega_2\Ad\,u_1\|<\eps$; by Kadison's transitivity 
one can make a small modification to $(u_t)$ to get the exact 
equality.

Obviously \ref{A10}$\Rightarrow$\ref{A7}$\Rightarrow$\ref{A8} 
among the properties defined above.

\begin{lem}\label{A105}
When $A$ is unital, Property \ref{A10} is equivalent to the one 
obtained by restricting the ambient algebra $B$ to a \cstar\ 
having $1_A$ as a unit. 
\end{lem} 
\begin{pf} 
Given $(\F,\eps)$ let $(\G,\delta)$ be the one obtained in this 
weaker property. Let $\G_1=\G\cup \{N1_A\}$, where $N$ will be 
specified later to be a large positive number, and suppose that 
$B$ is given such that $B\supset A$ and $1_A$ is not an identity 
for $B$. If $\pi$ is an irreducible representation of $B$ and 
$\om_1,\om_2$ are pure states of $B$ defined by the unit vectors 
$\xi_1,\xi_2\in \calH_{\pi}$ respectively such that  
$|\om_1(x)-\om_2(x)|<2^{-1}\eps^2\delta,\ x\in\G_1$, then we can 
choose $N$ so large  that $|\om_1(1_A)-\om_2(1_A)|$ is very small 
and either, $\om_1(1_A)$ and $\om_2(1_A)$ are smaller than 
$\eps^2$, or $|\phi_1(x)-\phi_2(x)|<\delta,\ x\in\G$ for the 
states $\phi_i=\om_i(1_A)^{-1}\om_i|1_AB1_A$. In any case we can 
choose a continuous path $(w_t)_{t\in[0,1]}$ in 
$\U((1-1_A)B(1-1_A))$ such that 
$\|(1-\pi(1_A))\xi_1-\pi(w_1)\xi_2\|$ is very small. In the former 
case where $\|\pi(1_A)\xi_i\|<\eps$, we set $v_t=1_A$ and in the 
latter case we choose, by the weaker property, a continuous path $ 
(v_t)_{t\in[0,1]}$ in $\U(1_AB1_A)$ such that 
$\|\pi(1_A)\xi_1-\pi(v_1)\xi_2\|$ is very small and 
$\|\Ad\,v_t(x)-x\|<\eps,\ x\in\F$. Setting $u_t=w_t+v_t$, it 
follows that $\|\xi_1-\pi(u_1)\xi_2\|$ is of order $\eps$ and 
$\|\Ad\,u_t(x)-x\|<\eps,\ x\in\F$. This completes the proof. 
\end{pf}

 \ncm{\calC}{{\mathcal{C}}}
\begin{prop}\label{A11}
If $\calC$ denotes the class of \cstars\ with Property \ref{A10}, 
then the following statements hold: 
 \begin{enumerate} 
 \item If $A\in\calC$ is non-unital, then $\tilde{A}\in\calC$, 
 where $\tilde{A}$ is the \cstar\ obtained by adjoining a unit to 
 $A$.
 \item If $A_1,A_2\in\calC$ are unital, then $A_1\oplus A_2\in\calC$. 
 \item If $A\in\calC$ and $e\in A$ is a projection, then 
 $eAe\in\calC$.
 \item If $A\in\calC$ and $I$ is an ideal of $A$, then the 
 quotient $A/I\in\calC$.
 \item If $(A_n)$ is an inductive system with $A_n\in\calC$, then 
 $\lim_n A_n\in\calC$.
 \end{enumerate}
\end{prop} 
\begin{pf}
To prove (1) we may assume that $1_B=1_{\tilde{A}}$ in Property 
\ref{A10}. Then this is obvious. 

To prove (2) let $\F$ be a finite subset of $A\equiv A_1\oplus 
A_2$ and $\eps>0$. Then there is a finite subset $\F_i\subset A_i$ 
for each $i=1,2$ such that if $\|\Ad\,u(x)-x\|<\eps,\ 
x\in\F_1\cup\F_2$, then $\|\Ad\,u(x)-x\|<\eps,\ x\in\F$, for any 
$u\in\U(B)$ with $B\supset A$ and $1_B=1_A$. Choose 
$(\G_i,\delta_i)$ for $(\F_i,\eps)$ as in Property \ref{A10}. We 
may suppose that $\delta_1=\delta_2\equiv \delta$ by multiplying 
each element of $\G_1$ by $\delta_2/\delta_1$. Let 
$\G=\G_1\cup\G_2\cup\{N p_1,N p_2\}$, where $p_i$ is the identity 
of $A_i$ in $B$ (and so $p_1+p_2=1$). Let $\pi$ be an irreducible 
representation of $B$ and let $\om_i$ be a pure state of $B$ 
defined by a unit vector $\xi_i\in\calH_{\pi}$ for $i=1,2$. 

By choosing $N$ sufficiently large, we may assume that if 
$|\om_1(x)-\om_2(x)|<\eps^2\delta/2, \ x\in\G$, then it follows 
that $|\om_1(p_i)-\om_2(p_i)|$ is very small for $i=1,2$, and both 
$\om_1(p_i)$ and $\om_2(p_i)$ are smaller than $\eps^2$ or 
 $$
 |\frac{\om_1(x)}{\om_1(p_i)}-\frac{\om_2(x)}{\om_2(p_i)}|<\delta,
 \ \ x\in \G_i.
 $$
In the former case we set $u_{it}=p_i$ and in the latter case we 
obtain a continuous path $(u_{it})$ in $\U(p_iBp_i)$ such that 
$u_{i0}=p_i$, and 
 \BE
 \|\Ad\,u_{it}(x)-x\|&<&\eps,\ \ x\in \F_i,\\ 
 \om_1(p_i)^{-1}\om_1|p_iBp_i&=&\om_2(p_i)^{-1}\om_2\Ad\,u_{i1}|p_iBp_i.
 \EE 
By multiplying $(u_{it})$ by a $\T$-valued continuous function, we 
may further suppose that 
$\pi(u_{i1})\xi_1=\|\pi(p_i)\xi_1\|\|\pi(p_i)\xi_2\|^{-1}\pi(p_i)\xi_2$, 
where $\|\pi(p_i)\xi_1\|\|\pi(p_i)\xi_2\|^{-1}$ is arbitrarily 
close to $1$. Set $u_t=u_{1t}+u_{2t}$;  then 
$\|\om_1-\om_2\Ad\,u_1\|$ is of the order of $\eps$ (since 
$\|\pi(u_1)\xi_1-\xi_2\|$ is at most of order $\eps$) and 
$\|\Ad\,u_t(x)-x\|<\eps,\ x\in\F_1\cup\F_2$. 

To prove (3) let $e$ be a projection in $A\in\calC$ and let $B$ be 
a \cstar\ with $B\supset eAe$ and $1_B=e$. Then there is a \cstar\ 
$D$ such that $eDe=B$ and  the diagram
 $$
 \begin{array}{ccc} 
 eAe & \subset & B\\  
 \cap& &\cap\\
 A&\subset& D
 \end{array}
 $$
is commutative. (To show this we may suppose that $A$ acts on a 
Hilbert space $\calH$ and that $B$ acts on $e\calH$ with $B\supset 
eAe$. We set $D$ to be the \cstar\ generated by $B$ and $A$.) Let 
$\F$ be a finite subset of $eAe$ and $\eps>0$. Let 
$\F_1=\F\cup\{Ne\}$ with $N$ sufficiently large. We choose 
$(\G,\delta)$ for $(\F_1,\eps)$ from Property \ref{A10} of $A$ and 
set $\G_1=\{exe\ | \ x\in\G\}$. Let $\om_1,\om_2$ be pure states 
of $B$ such that $\om_1\sim\om_2$ and 
$|\om_1(x)-\om_2(x)|<\delta,\ x\in\G_1$. Extend $\om_i$ to a pure 
state of $D$; since the extension is unique, we denote it by the 
same symbol. Since $\om_i(x)=\om_i(exe)$, we obtain, from Property 
\ref{A10}, a continuous path $(u_t)$ in $\U(D)$ such that $u_0=1$, 
$\om_1=\om_2\Ad\,u_1$, and $\|\Ad\,u_t(x)-x\|<\eps,\ x\in\F_1$. 
Since $\|[e,u_t]\|$ is very small, we can obtain a continuous path 
in $\U(eAe)$ from $(eu_te)$ which satisfies the required 
properties. 

To prove (4) let $B$ a \cstar\ with $B\supset A/I$. Let $D=B\oplus 
A$ and embed $A$ into $D$ by $x\mapsto (x+I,x)$. Let $\F$ be a 
finite subset of $A/I$ and $\eps>0$. For each $x\in\F$ choose an 
$x_1\in A$ such that $x_1+I=x$ and denote by $\F_1$ the subset 
consisting of such $x_1$ with $x\in\F$. For $(\F_1,\eps)$ we 
obtain $(\G_1,\delta)$ as in Property \ref{A10}. Denoting by $p$ 
the projection of $D$ onto $B$, we set $\G=\{p(x)\ |\ x\in\G_1\}$. 
Then it is easy to check that $(\G,\delta)$ satisfies the requires 
properties for $(\F,\eps)$. 

To show (5) note, by (4), that $(A_n)$ can be regarded as an 
increasing sequence in $\calC$. Then for any finite subset $\F$ of 
$\ol{\cup_nA_n}$ we find $A_n$ which almost contains $\F$. Hence 
this is immediate. 
\end{pf}

\section{Homogeneous $C^*$-algebras}

In this section we will show that the \cstars\ of the form 
$C\otimes M_n$ have Property \ref{A10}, where $C$ is a unital 
abelian \cstar. Then it will follow by \ref{A11} that any 
approximately homogeneous \cstar\ has Property \ref{A10}. 
Furthermore we will attempt to prove that some subhomogeneous 
\cstars\ have Property \ref{A10}. 

The following is well-known.

\begin{lem}\label{C1}
For any $\eps>0$ there is a $\delta>0$ satisfying: If $A$ is a 
\cstar\ and $A_{+1}=\{x\in A\ |\ x\geq0,\,\|x\|\leq1\}$, then that 
$\|x-y\|<\delta$ implies that $\|x^{1/2}-y^{1/2}\|<\eps$ for any 
$x,y\in A_{+1}$. \end{lem} 
\begin{pf}
For any $\eps>0$ there is a real-valued polynomial $p(t)$ with 
$p(0)=0$ such that $|t^{1/2}-p(t)|<\eps,\ t\in[0,1]$. If 
$p(t)=\sum_{i=1}^na_it^i$, set $C=\sum_i|a_i|i$. Then for $x,y\in 
A_{+1}$ we have that $\|p(x)-p(y)\|\leq C\|x-y\|$. Hence for 
$x,y\in A_{+1}$ the condition that $\|x-y\|<\eps/C$ implies that 
$\|x^{1/2}-y^{1/2}\|<3\eps$. \end{pf}   

\begin{lem}\label{C2}
For any $\eps>0$ and $n\in\N$ there exists a $\delta>0$ 
satisfying: If $\xi_1,\xi_2,\ldots,\xi_n$ are vectors in a Hilbert 
space $\calH$ with $\dim(\calH)\geq n$ and $c=(c_{ij})$ is an 
$n\times n$ matrix such that $\sum_{i=1}^n\|\xi_i\|^2\leq1$, 
$c\geq0$, and 
 $$
 |c_{ij}-\lan\xi_i,\xi_j\ran|<\delta,\ \ i,j=1,2,\ldots,n,
 $$
then there exist vectors $\eta_1,\eta_2,\ldots,\eta_n$ in $\calH$ 
such that
 \BE
 \lan\eta_i,\eta_j\ran&=&c_{ij},\ \ i,j=1,2,\ldots,n,\\
 \|\xi_i-\eta_i\|&<&\eps,\ \ i=1,2,\ldots,n. 
 \EE
\end{lem}  
\begin{pf}
Let $d_{ij}=\lan\xi_i,\xi_j\ran$. Then the $n\times n$ matrix 
$d=(d_{ij})$ satisfies that $d\geq0$ and $\|d\|\leq1$. The 
condition that $|c_{ij}-d_{ij}|<\delta$ for all $i,j$ implies that 
$\|c-d\|\leq n\delta$.

If $d$ is strictly positive, then define 
$\eta_i=\sum_{k=1}^n(c^{1/2}d^{-1/2})_{ik}\xi_k$. Then by 
computation
 $$\lan\eta_i,\eta_j\ran=\sum_{k,\ell}(c^{1/2}d^{-1/2})_{ik} 
 \ol{(c^{1/2}d^{-1/2})}_{j\ell}d_{k\ell}=c_{ij}
 $$
and
 $$
 \|\eta_i-\xi_i\|^2=(c^{1/2}-d^{1/2})^2_{ii}\leq \|c^{1/2}-d^{1/2}\|^2.
 $$
By the previous lemma one can choose $\delta>0$ so small that 
$\|c^{1/2}-d^{1/2}\|<\eps$ follows from $\|c-d\|\leq n\delta$ and 
$\|d\|\leq1$.

In general let $\calL$ be a subspace of $\calH$ such that 
$\calL\ni\xi_i$ for all $i$ and $\dim\calL=n$. Then there is a 
sequence $(\xi_{k1},\xi_{k2},\ldots,\xi_{kn})$ of bases of $\calL$ 
such that $\|\xi_{ki}-\xi_i\|\ra0$ as $k\ra\infty$. For each 
$(\xi_{k1},\xi_{k2},\ldots,\xi_{kn})$ one can apply the previous 
argument to produce vectors $\eta_{k1},\eta_{k2},\ldots,\eta_{kn}$ 
in $\calL$. By using the compactness argument for the limit of $k$ 
to infinity one obtains vectors $\eta_1,\eta_2,\ldots,\eta_n$ in 
$\calL$ such that $\lan\eta_i,\eta_j\ran=c_{ij}$ and 
$\|\eta_i-\xi_i\|^2=(c^{1/2}-d^{1/2})^2_{ii}$. This completes the 
proof. \end{pf} 

\begin{lem}\label{C3}
For any $\eps>0$ and $n\in\N$ there exists a $\delta>0$ 
satisfying: Given sequences $(\xi_1,\xi_2,\ldots,\xi_n)$ and 
$(\eta_1,\eta_2,\ldots,\eta_n)$ of vectors in a Hilbert space 
$\calH$ such that $\sum_i\|\xi_i\|^2\leq1$, 
$\sum_i\|\eta_i\|^2\leq1$, and
 $$
 |\lan\xi_i,\xi_j\ran-\lan\eta_i,\eta_j\ran|<\delta,\ \ i,j=1,2,\ldots,n,
 $$
there is a unitary $U$ on $\calH$ such that
 $$
 \|U\xi_i-\eta_i\|<\epsilon,\ \ i=1,2,\ldots,n.
 $$
\end{lem}
\begin{pf}
If $\dim\calH\geq n$, then choose $\delta>0$ as in the previous 
lemma and find a sequence $(\zeta_1,\zeta_2,\ldots,\zeta_n)$ in 
$\calH$ such that $\lan\zeta_i,\zeta_j\ran=\lan\xi_i,\xi_j\ran$ 
and $\|\zeta_i-\eta_i\|<\epsilon$. Then we can find a unitary $U$ 
on $\calH$ such that $U\xi_i=\zeta_i$ for $i=1,2,\ldots,n$, which 
satisfies the required properties. 

If $m=\dim\calH<n$, we first choose $m$ vectors from  
$(\xi_1,\xi_2,\ldots,\xi_n)$. Pick up $\xi_{i_1}$ with 
$\|\xi_{i_1}\|=\max_j\|\xi_j\|$. Let $P_1$ be the projection onto 
the subspace spanned by $\xi_{i_1}$ and pick up $\xi_{i_2}$ with 
$\|(1-P_1)\xi_{i_2}\|=\max_j\|(1-P_1)\xi_j\|$. Let $P_2$ be the 
projection onto the subspace spanned by $\xi_{i_1},\xi_{i_2}$ and 
pick up $\xi_{i_3}$ with 
$\|(1-P_2)\xi_{i_3}\|=\max_j\|(1-P_2)\xi_j\|$. Repeating this 
process we obtain $\xi_{i_1},\xi_{i_2},\ldots,\xi_{i_m}$, which we 
assume are all different as we may. If $j\not\in 
I=\{i_1,i_2,\ldots,i_m\}$, then there are $c^k_j\in\C$ for $k\in 
I$ such that
 $$
 \xi_j=\sum_{k\in I}c_j^k\xi_k.
 $$
From the construction of $I$ we can assume that $|c_j^k|\leq1$. 
(If $\xi_k,\,k\in I$ are linearly independent, then of course 
$|c_j^k|\leq1$ automatically.)  If 
$|\lan\xi_i,\xi_j\ran-\lan\eta_i,\eta_j\ran|<\delta$ for all 
$i,j$, then $\|\eta_j-\sum_{k\in I}c_j^k\eta_k\|^2\leq 
(m+1)^2\delta$. Thus if we define a unitary $U$ on $\calH$ by 
requiring that $\|U\xi_k-\eta_k\|<\eps,\ k\in I$, we obtain that 
for $j\not\in I$,
 $$
 \|U\xi_j-\eta_j\|\leq \sum_{k\in I}|c_j^k|\|U\xi_k-\eta_k\|
 +\|\eta_j-\sum_{k\in I}c_j^k\eta_k\|\leq m\eps+(m+1)^2\delta.
 $$
This completes the proof. \end{pf}
                                     
\begin{lem}\label{C4}
Any matrix algebra $M_n$ has Property \ref{A10}. 
\end{lem} 
\begin{pf}
Let $B$ be a \cstar\ containing $A=M_n$ with $1_B=1_A$. Let 
$(e_{ij})$ be the set of matrix units of $A=M_n$. Since $A$ is 
spanned by the matrix units, we may take $\{e_{ij}\ |\ 
i,j=1,2,\ldots,n\}$ for both $\F$ and $\G$ in Property \ref{A10}.  
Let $\eps>0$ and let $\delta>0$ be the $\delta$ for $\eps 
n^{-1/2}$ in place of $\eps$ in Lemma \ref{C3}.

Let $\pi$ be an irreducible representation of $B$ on $\calH$ and 
let $\xi,\eta$ be unit vectors in $\calH$ such that
 $$
 |\lan\pi(e_{ij})\xi,\xi\ran-\lan\pi(e_{ij})\eta,\eta\ran|<\delta
 $$
for all $i,j$. We apply Lemma \ref{C3} to the sequences 
$(\pi(e_{1j})\xi)$ and $(\pi(e_{1j})\eta)$ and obtain a unitary 
$U$ on the subspace $\calL$ spanned by $(\pi(e_{1j})\xi)$ and 
$(\pi(e_{1j})\eta)$ such that 
$\|U\pi(e_{1j})\xi-\pi(e_{1j})\eta\|<\eps n^{-1/2}$. Since 
$U=e^{iH}$ with $\|H\|\leq\pi$ we choose a self-adjoint $h=h^*\in 
e_{11}Be_{11}$ such that  $\|h\|\leq\pi$ and $\pi(h)|\calL=H$. 
Define a unitary $u_t$ for each $t\in[0,1]$ by 
 $$
 u_t=\sum_ie_{i1}e^{\sqrt{-1}th}e_{1i}.
 $$
Then it follows that $u_t\in B\cap A'$. By computation
 \BE
 \|\pi(u_1)\xi-\eta\|^2&=&\|\sum_i\{\pi(e_{i1})U\pi(e_{1i})\xi-\pi(e_{ii})\eta\}
           \|^2\\ 
     &=&\sum_i\|U\pi(e_{1i})\xi-\pi(e_{1i})\eta\|^2\\     
    &<& \eps^2.
  \EE
Thus $\|\pi(u_1)\xi-\eta\|<\eps$. Since $u_t\in B\cap A'$, this 
completes the proof. \end{pf}
           
\begin{prop}\label{C5}
If $A$ is a unital \cstar\ satisfying Property \ref{A10} and $C$ 
is a unital commutative \cstar, then $A\otimes C$ satisfies 
Property \ref{A10} \end{prop}

Since any unital commutative \cstar\ is an inductive limit of 
quotients  of $C(\T^n)$, we may assume that $C=C(\T)$ in the above 
proposition by \ref{A11}. This follows essentially from the 
following lemma. 

\begin{lem}\label{C6}
For any small $\eps>0$ and $\eps'>0$ there exist a $\delta>0$ and 
$n\in\N$ satisfying: Let $U$ be a unitary operator on a Hilbert 
space $\calH$ with $E$ its spectral measure and let $\xi,\eta$ be 
two unit vectors in $\calH$. If 
 $$|\lan U^k\xi,\xi\ran-\lan U^k\eta,\eta\ran| <\delta
 $$
for $k=0,\pm1,\ldots, \pm n$, there exists a sequence 
$(t_1,t_2,\ldots,t_m)$ of points in $\T$ such that $t_{i-1}<t_i$ 
in the cyclic order with $t_0=t_m$, 
$\eps/2<\dist(t_{i-1},t_i)<3\eps/2$ and 
 \BE
 && |\|E(t_{i-1},t_i)\xi\|^2-\|E(t_{i-1},t_i)\eta\|^2|<2\eps',\\
 &&\|E(t_{i}-\gamma/2,t_i+\gamma/2)\xi\|^2<\eps',\ \ 
 \|E(t_{i}-\gamma/2,t_i+\gamma/2)\eta\|^2<\eps',
 \EE
where $\gamma=\eps\eps'/4$ and the total length of $\T$ is 
normalized to be $1$. 
\end{lem} 
\begin{pf}
Let $\gamma=\eps\eps'/4$. For any interval $I$ of $\T=\R/\Z$ of 
length $\eps/2$, there is a $t\in I$ such that both 
 $\|E(t+(-\gamma/2,\gamma/2])\xi\|^2
 $
and  
 $\|E(t+(-\gamma/2,\gamma/2])\eta\|^2
 $ 
are smaller than $\eps'$. Otherwise
 $$\lan E(t+(-\gamma/2,\gamma/2])\xi,\xi\ran+\lan 
 E(t+(-\gamma/2,\gamma/2])\eta,\eta\ran\geq \eps'
  $$
for any $t\in I$, which implies that 
 $$
 \lan E(J)\xi,\xi\ran+\lan E(J)\eta,\eta\ran>2,
 $$                             
for an interval $J$ including $I$, a contradiction. Thus there is 
a sequence $(t_1,t_2,\ldots,t_{m})$ of points in $\T$ such that 
$t_{i-1}<t_i$ in the cyclic order and  
$\eps/2<\dist(t_{i-1},t_i)<3\eps/2$ with $t_0=t_m$, and 
$\|E(J_i)\xi\|^2<\eps'$ and  $\|E(J_i)\eta\|^2<\eps'$ for 
$J_i=t_i+(-\gamma/2,\gamma/2]$. 

Let $f$ be a $C^{\infty}$-function on $\T$ such that $0\leq f\leq 
1$, $\supp f\subset [0,1/4]$, and $f=1$ on 
$[\gamma/2,1/4-\gamma/2]$. If $\xi,\eta$ satisfy the condition in 
the statement  for a sufficiently small $\delta$ and for a 
sufficiently large $n$, it follows that for any product $g$ of two 
translates of $f$, 
 $$
 |\lan g(U)\xi,\xi\ran-\lan g(U)\eta,\eta\ran|<\eps'.
 $$
We construct a function $f_i$, as the product of two translates of 
$f$, such that $\supp f_i\subset [t_{i-1},t_i]$ and $f_i=1$ on 
$[t_{i-1}+\gamma/2,t_i-\gamma/2]$. Since $|\lan 
f_i(U)\xi,\xi\ran-\lan f_i(U)\eta,\eta\ran|<\eps'$, and 
$\|E(J_i)\xi\|^2<\eps'$ etc., we have that for any $i$, 
 $$
 |\|E(t_{i-1},t_i)\xi\|^2-\| E(t_{i-1},t_i)\eta\|^2|<2\eps'.
 $$                                                 
This completes the proof. 
\end{pf} 
         
\medskip \noindent  
{\em Proof of Proposition \ref{C5}}\ \ By the following lemma, for 
any finite subset $\F$ of $A$ and $\eps>0$ we have a finite subset 
$\G_A$ of $A$ and $\delta'>0$ satisfying: Let $B$ be a non-unital 
\cstar\ such that $A$ is a unital $C^*$-subalgebra of the 
multiplier algebra $M(B)$ of $B$. If two pure states $\om_1,\om_2$ 
of $B$ satisfy that $\om_1\sim\om_2$ and 
$|\om_1(x)-\om_2(x)|<\delta',\ x\in \G_A$, where $\om_i$ also 
denotes the natural extension of $\om_i$ to a state on $M(B)$, 
then there is a  continuous path $(u_{t})$ in $\U(B)$ such that 
$u_{t}-1\in B$,  $u_{0}=1$, $\om_1=\om_2\Ad\,u_{1}$, and 
 $$
 \|\Ad\,u_{t}(x)-x\|<\eps,\ \ x\in\F.
 $$
We may assume, by replacing $\delta$ by a smaller one, that 
$\|x\|\leq1$ for $x\in\G_A$ and also that $1\in\G_A$.                                           

Suppose that $B\supset A\otimes C(\T)$ with the common unit. Given 
$\eps>0$  and $\eps'=\eps^3\delta'/4$ we choose $\delta,n$ as in 
the previous lemma. 

Let $\G=\{xz^k\ |\ x\in\G_A,\ |k|\leq n\}$, where $z$ is the 
canonical unitary of $C(\T)\simeq1\otimes C(\T)\subset A\otimes 
C(\T)$. Let $\pi$ be an irreducible representation of $B$ and 
$\om_1,\om_2$ be the pure states of $B$ defined through  unit 
vectors $\xi,\eta$ respectively. Suppose that 
 $$
 |\om_1(x)-\om_2(x)|<\min(\delta,\eps'), \ \ x\in\G.
 $$  

  For $U=\pi(z)$ we have that 
 $$
  |\lan U^k\xi,\xi\ran-\lan U^k\eta,\eta\ran| <\delta, \ \ |k|\leq n.
 $$
Then we obtain a sequence $(t_1,t_2,\ldots,t_m)$ in $\T$ as in 
\ref{C6}. Let $B_i$ be the hereditary $C^*$-subalgebra of $B$ 
generated by $E(t_{i-1},t_i)$, where $E$ is the spectral measure 
of $z$. There is a natural unital embedding of $A$ into the 
multiplier algebra $M(B_i)$. Let $f_i\in C(\T)$ be a function 
supported on $[t_{i-1},t_i]$ as given in the proof of the previous 
lemma. Then from the proof of \ref{C6} we may assume that 
 $$
 |\om_1(xf_i(z))-\om_2(xf_i(z))|<\eps',\ \ x\in \G_A,
 $$  
which then implies that 
 $$
 |\lan xE(t_{i-1},t_i)\xi,\xi\ran-
 \lan xE(t_{i-1},t_i)\eta,\eta\ran|<2\eps', \ \ 
 x\in\G_A.
 $$

We define constants $c_i,d_i$ by 
 $$d_i\|E(t_{i-1},t_i)\xi\|=c_i\| 
 E(t_{i-1},t_i)\eta\|= \min(\|E(t_{i-1},t_i)\xi\|,\| 
 E(t_{i-1},t_i)\eta\|),$$ 
where  $E$ is regarded as the spectral measure of $U=\pi(z)$.  If 
$$d_i\|E(t_{i-1},t_i)\xi\|=c_i\| 
E(t_{i-1},t_i)\eta\|\leq\eps^{3/2},
$$ then we set 
$u_{it}=1\in\U(B_i)$; otherwise for  the states on $B_i$: 
 $\varphi_1=\|E(t_{i-1},t_i)\xi\|^{-2}\om_1|B_i$ and 
 $\varphi_2=\|E(t_{i-1},t_i)\eta\|^{-2}\om_2|B_i$,  it follows 
that 
 $$|\varphi_1(x)-\varphi_2(x)|<4\eps'\eps^{-3}=\delta',\ \ x\in\G_A.
 $$           
Here we have used the fact that 
$|1-\|E(t_{i-1},t_i)\xi\|^{-2}\|E(t_{i-1},t_i)\eta\|^2|<2\eps'\eps^{-3}$
and $\|x\|\leq1$ for $x\in\G_A$. Hence we obtain a continuous path 
$(u_{it})$ in $\U(B_i)$ such that $u_{it}-1\in B_i$, 
$\|\Ad\,u_{it}(x)-x\|<\eps,\, x\in\F$, and  $ 
\pi(u_{i1})d_iE(t_{i-1},t_i)\xi=c_iE(t_{i-1},t_i)\eta$. Then the 
sum $(u_t)$ of $(u_{it})$ over $i$ defines a continuous path in 
$\U(B)$ (or regarding each $u_{it}$ as in $\U(B)$ take the product 
of these $u_{it}$ over $i$). Then, since 
$\dist(t_{i-1},t_i)<3\eps/2$, we have that 
$\|u_tz-zu_t\|<3\pi\eps$ and  also that 
 $$
 \|\Ad\,u_{t}(x)-x\|<\eps,\ \ x\in\F.
 $$
Let $S$ be the set of $i$ with $d_i\|E(t_{i-1},t_i)\xi\|> 
\eps^{3/2}$. Then  
 $$
 \sum_{i\not\in S}\|d_iE(t_{i-1},t_i)\xi\|^2<2\eps^2,
 $$
and hence
 $$
 \sum_{i\in S}\|d_iE(t_{i-1},t_i)\xi\|^2>1-3\eps^2,
 $$                                                
because 
$\sum_i\|E(\{t_{i}\})\xi\|^2<2\eps'\eps^{-1}<\eps^2\delta'/2$, 
which can be assumed to be smaller than $\eps^2$. Since 
 $
 \pi(u_1)d_iE(t_{i-1},t_i)\xi=c_iE(t_{i-1},t_i)\eta,
 $
for $i\in S$, we have that
 $$
 \pi(u_1)\sum_{i\in S}d_iE(t_{i-1},t_i)\xi=\sum_{i\in 
 S}c_iE(t_{i-1},t_i)\eta.
 $$
This implies that 
 $$ \|\pi(u_1)\xi-\eta \|<2\sqrt{3}\eps.
 $$
This completes the proof.  \hfill $\square$ 

\medskip 
 
\begin{lem} \label{mult}
If $A$ is a unital \cstar\ satisfying Property \ref{A10}. Then for 
any finite subset $\F$ of $A$ and $\eps>0$ there is a finite 
subset $\G$ of $A$ and $\delta>0$ satisfying: If $B$ is a 
non-unital \cstar\ such that  $A\subset M(B)$ with common unit and 
$\pi$ is an irreducible representation of $B$, and $\xi,\eta$ are 
unit vectors in $\calH_{\pi}$ such that if  
 $$ |\lan\pi(x)\xi,\xi\ran-\lan\pi(x)\eta,\eta\ran|<\delta,\ \ x\in\G,
 $$                                  
where $\pi$ also denotes the natural extension of $\pi$ to a 
representation of $M(B)$, then there is a continuous path 
$(u_t)_{t\in[0,1]}$ in $\U(B)$ such that $u_0=1$, 
$\eta=\pi(u_1)\xi$, $u_t-1\in B$, and 
 $$ \|\Ad\,u_t(x)-x\|<\eps,\ \ x\in\F,\ t\in [0,1].
 $$       
\end{lem}
\begin{pf}
In the situation as above, for the choice of $(\G,\delta)$ made 
for Property \ref{A10}, we have a continuous path $(u_t)$ in 
$\U(M(B))$ such that $u_0=1$, $\eta\in\C\pi(u_1)\xi$ and 
 $$ \|\Ad\,u_t(x)-x\|<\eps,\ \ x\in\F,\ t\in [0,1].
 $$ 
By multiplying $(u_t)$ by a suitable continuous function, we may 
assume that $\eta=\pi(u_1)\xi$.

Let $(e_{\iota})$ be an approximate identity of $B$. Since 
$[e_{\iota},x]\in B$ converges to zero in the $\sigma(B,B^*)$ 
topology for any $x\in M(B)$, we can assume, by taking a net from 
the convex combinations of $e_{\iota}$'s if necessary, that 
$[e_{\iota},x]$ converges to zero in norm for $x\in M(B)$; in 
particular $[e_{\iota},u_t]$ converges to zero as well as 
$[e_{\iota},x],\ x\in\F$. Thus we find a sequence $(e_n)$ in $B$ 
such that $0\leq e_n\leq 1,\ \pi(e_n)\xi\ra\xi,\  [e_n,u_t]\ra0$,  
$ [e_n,x]\ra0,\ x\in\F$, and $e_ne_m-e_m\ra0$ as $n\ra\infty$ for 
each $m$. Let $\mu\in (0,1)$ be so close to 1 that $\|u_t-u_{\mu 
t}\|<\eps$ for $t\in[0,1]$ and let $k\in\N$ be such that 
$\|u_{\mu^kt}-1\|<\eps$. For a subsequence $(n_i)$ set 
 $$
 z_t=u_te_{n_1}+u_{\mu t}p_1+u_{\mu^2 t}p_2+\cdots +u_{\mu^k 
 t}p_k+1-e_{n_{k+1}},
 $$
where $p_i=e_{n_{i+1}}-e_{n_i}$. By assuming that 
$e_{n_{i+1}}e_{n_i}\approx e_{n_i}$, we have that $p_ip_j\approx0$ 
if $|i-j|>1$. We can also assume that $z_t$ almost commutes with 
$e_{n_i}$ and $x\in\F$. Hence $z_tp_i\approx u_{\mu^it}p_i$, 
$z_te_{n_1}\approx u_te_{n_1}$, and $z_t(1-e_{n_{k+1}})\approx 
1-e_{n_{k+1}}$, all up to $\eps$. Since $1=e_{n_1}+\sum_{i=1}^k 
p_i+ 1-e_{n_{k+1}}$ and $p_iz_t^*z_tp_j$ is arbitrarily close to 
zero if $|i-j|>1$ for $i,j=0,1,\ldots,k+1$ with 
$p_0=e_{n_1},\,p_{k+1}=1-e_{n_{k+1}}$, we can conclude that 
$z_t^*z_t\approx 1$ up to the order of $\eps$. (For example 
$\|z_t^*z_t-1\|=\|\sum_{i=0}^{k+1}(z_t^*z_t-1)p_i\|\leq 
\|\sum(z_t^*z_t-1)p_{2i}\|+\|\sum(z_t^*z_t-1)p_{2i+1}\|\lesssim2\max_i\|(z_t^*z_t-1)p_i\|$.) 
Hence $z_t$ is close to a unitary in $B+1$ up to the order of 
$\eps$. By assuming that $\pi(e_{n_1})\xi\approx\xi$, we have that 
$\pi(z_1)\xi\approx\eta$. In this way we can construct a 
continuous path $(v_t)$ in $\U(B)$ in a small neighborhood of 
$(z_t)$ (of order $\eps$) such that $\pi(v_1)\xi\approx\eta$, 
$v_t-1\in B$, and $\|[v_t,x]\|\approx0$ up to the order of $\eps$ 
for $x\in\F$. This completes the proof. 
\end{pf}

For $n=1,2,\ldots$ let $D_n$ denote the  dimension drop \cstar: 
 $$
 D_n=\{x\in C([0,1];M_n)\ |\ x(0),x(1)\in\C1_n\}.
 $$
 
\begin{lem}\label{drop}
For any $k,n\in\N$, $D_n\otimes M_k$ has Property \ref{A10}. 
\end{lem} 
\begin{pf}
For any finite subset $\F$ of $D_n\otimes M_k$ and $\eps>0$, we 
find an interval $[a,b]\subset (0,1)$ such that if $x\in\F$, then 
$\|x(t)-x(0)\|<\eps$ for any  $t\in[0,a]$ and $\|x(t)-x(1)\|<\eps$ 
for any $t\in [b,1]$. Hence we may suppose that if $x\in\F$, then 
$x(t)=x(0)\in 1_n\otimes M_k$ for $t\in[0,a]$ and $x(t)=x(1)\in 
1_n\otimes M_k$ for $t\in [b,1]$. Replacing $a$ by a smaller one 
and $b$ by a larger one, we may further suppose that if $x\in\F$, 
then $x(t)=x(0)$ around $t=a$ and $x(t)=x(1)$ around $t=b$. 
Denoting by $\F|[a,b]$ the subset of $C[a,b]\otimes M_n\otimes 
M_k$ obtained by restricting $\F$ to $[a,b]$, we choose 
$(\G_2,\delta_2)$ for $(\F|[a,b],\eps)$ as in \ref{mult}. 
Similarly we choose $(\G_0,\delta_0)$ for $(\F|[0,a],\eps)$ with 
$\F|[0,a]$ as a subset of $1\otimes M_k$ and $(\G_1,\delta_1)$ for 
$(\F|[b,1],\eps)$ with $\F|[b,1]$ as a subset of $1\otimes M_k$ by 
using \ref{mult}. We may assume that 
$\delta_0=\delta_1=\delta_2\equiv\delta$ and $\|x\|\leq1$ for 
$x\in\G_0\cup\G_1\cup\G_2$. Furthermore we may assume that if 
$x\in\G_2$, then $x(t)=x(a)$ around $t=a$ in $[a,b]$ and 
$x(t)=x(b)$ around $t=b$ in $[a,b]$. For $x\in\G_2$, we set 
 $$
 \ol{x}(t)=\left\{\begin{array}{ll} x(b)&t>b\\ x(t)&a\leq t\leq b\\
                                    x(a)&t<a \end{array}\right.
 $$
and define $\ol{\G}_2=\{\ol{x}\ |\ x\in\G_2\}\cup \{1\}\subset 
C[0,1]\otimes M_n\otimes M_k$. If $0<a'<a$ and $b<b'<1$, we can 
choose $(\ol{\G}|[a',b'],\delta)$ as the $(\G,\delta)$ for 
$(\F|[a',b'],\eps)$ in \ref{mult}; because there is an isomorphism 
of $C[a,b]\otimes M_n\otimes M_k$ onto $C[a',b'] \otimes 
M_n\otimes M_k$ which sends $\F|[a,b]$ onto $\F|[a',b']$ and 
$\G_2=\ol{\G}_2|[a,b]$ onto $\ol{\G}_2|[a',b']$ respectively. 

Let $N\in\N$ be such that $N\delta\eps^2>16$, and let $\eps'>0$ be 
such that $N\eps'<\min(a,1-b)$. Let 
 $$
 I_k=(a-(k+1)\eps',a-k\eps'],\ \ J_k=(b+k\eps',b+(k+1)\eps']
 $$
for $k=0,1,\ldots,N-1$. Let $f_k\in C[0,1]$ be such that $0\leq 
f_k\leq 1$, $\supp f_k\subset (a-(k+1)\eps',b+(k+1)\eps']$, and 
$f_k(t)=1$ for $t\in [a-k\eps',b+k\eps']$ and let 
$f_{0k}=(1-f_k)\chi_{[0,a]}\in C[0,1]$ and 
$f_{1k}=(1-f_k)\chi_{[b,1]}\in C[0,1]$. Finally we define a subset 
$\G$ of $D_n\otimes M_k$ as the union of
 $
 \{f_kx\ |\ x\in\ol{\G}_2,\,k=0,1,\ldots,N-1\}$,
 $ \{f_{0k}x\ |\ x\in \G_0\ {\rm or}\ x=1;\, k=0,\ldots,N-1\}$, and 
 $\{f_{1k}x\ |\ x\in \G_1\ {\rm or}\ x=1;\,k=0,\ldots,N-1\}$.

Let $B$ be a \cstar\ containing $D_n\otimes M_k$ as a unital 
$C^*$-subalgebra and let $\pi$ be an irreducible representation of 
$B$. Let $\xi,\eta$ be unit vectors in $\calH_{\pi}$ such that 
 $$
 |\lan\pi(x)\xi,\xi\ran-\lan\pi(x)\eta,\eta\ran|<\delta\eps^2/4,\ \ x\in \G.
 $$
We first choose $k$ such that 
 $$
 s_k\equiv \|\chi_{I_k}\xi\|^2+\|\chi_{I_k}\eta\|^2+
            \| \chi_{J_k}\xi\|^2+\|\chi_{J_k}\eta\|^2
             <\delta\eps^2/4,
 $$
where $\chi_{I_k}$ is the characteristic function of $I_k$ 
regarded as an element of $(D_n\otimes M_k)^{**}\subset B^{**}$, 
which is also identified with $\pi^{**}(\chi_{I_k})$ etc. 
(Otherwise $s_k\geq \delta\eps^2/4$ for all $k$ leads us to a 
contradiction because $N\delta\eps^2>16$.)

Let $a'=a-(k+1)\eps'$ and $b'=b+(k+1)\eps'$, and let 
$e=\chi_{(a',b')}\in B^{**}$ and $B_2$ the hereditary 
$C^*$-subalgebra of $B$ corresponding to $e$. Since 
 $
 |\lan\pi(f_kx)\xi,\xi\ran-\lan\pi(f_kx)\eta,\eta\ran|<\delta\eps^2/4,\ 
 \ x\in \ol{\G}_2,
 $
we have that 
 $$
 |\lan ex\xi,\xi\ran-\lan ex\eta,\eta\ran|<\delta\eps^2/2,\ 
 \ x\in\ol{\G}_2,
 $$              
where we have omitted $\pi$ and $\pi^{**}$.  If $\|e\xi\|>\eps$ 
(and $\eps,\delta$ are sufficiently small), then it follows that
 $$
 |\lan\pi_2(x)\xi_2,\xi_2\ran-\lan\pi_2(x)\eta_2,\eta_2\ran|<\delta,
 \ \ x\in \ol{\G}_2|[a',b'], 
 $$                         
where $\pi_2$ is the irreducible representation of $B_2$ on 
$e\calH_{\pi}$  obtained by restricting $\pi$,  and 
$\xi_2=\|e\xi\|^{-1}e\xi$ and $\eta_2=\|e\eta\|^{-1}e\eta$.  Note 
that $C[a',b']\otimes M_n\otimes M_k$ is regarded as a subalgebra 
of the multiplier algebra of $B_2$.  Hence by \ref{mult} we can 
find a continuous path $(v_t)$ in $\U(B_2)$ such that $v_0=1$, 
$v_t-1\in B_2$, $\pi_2(v_1)\xi_2= \eta_2$, and 
$\|\Ad\,v_t(x)-x\|<\eps,\ x\in \F|[a',b']$.  If $\|e\xi\|\leq 
\eps$, we set $v_t=1$. In any case we have that 
$\pi_2(v_1)\xi_2\approx\eta_2$ up to $\eps$. 

Let $B_0$ be the hereditary $C^*$-subalgebra of $B$ corresponding 
to $e_0=\chi_{[0,a')}$. If $\|e_0\xi\|>\eps$, it follows that 
 $$
 |\lan\pi_0(x)\xi_0,\xi_0\ran-\lan\pi_0(x)\eta_0,\eta_0\ran|<\delta,
 \ \ x\in \G_0,
 $$
where $\pi_0$ is the irreducible representation of $B_0$ on 
$e_0\calH_{\pi}$ obtained by restricting $\pi$, and 
$\xi_0=\|e_0\xi\|^{-1}e_0\xi$ and 
$\eta_0=\|e_0\eta\|^{-1}e_0\eta$. Hence we obtain a continuous 
path $(v_{0t})$ in $\U(B_0)$ such that $v_{0,0}=1$, $v_{0t}-1\in 
B_0$, $\pi_0(v_{0,1})\xi_0= \eta_0$, and 
$\|\Ad\,v_{0t}(x)-x\|<\eps,\ x\in \F|[0,a']$. If 
$\|e_0\xi\|\leq\eps$, then we set $v_{0t}=1$. As before in any 
case we have that $\pi_0(v_{0,1})\xi_0\approx\eta_0$ up to $\eps$. 

In a similar way  we obtain a continuous path $(v_{1t})$ in 
$\U(B_1)$ with $B_1$ the hereditary $C^*$-subalgebra of $B$ 
corresponding to $\chi_{(b',1]}$ with similar properties to the 
above. Since 
$\|(e+e_0+e_1)\xi\|^2=\|e\xi\|^2+\|e_0\xi\|^2+\|e_1\xi\|^2\approx1$ 
up to $\eps^2$ , we get the desired path in $\U(B)$ by summing 
these paths $(v_t),(v_{0t})$, and $(v_{1t})$. 
\end{pf}

We recall that $A$ is called an AH algebra (or an approximately 
homogeneous \cstar) if $A$ is an inductive limit of \cstars\ of 
the form $\oplus_{i=1}^ke_i(C_i\otimes M_{n_i})e_i$, where $C_i$ 
is a unital commutative \cstar\ and $e_i$ is a projection of 
$C_i\otimes M_{n_i}$. The class of AH algebras include UHF 
algebras and AF algebras. 

More generally, a \cstar\ $A$ is called an ASH algebra (or an 
approximately sub-homogeneous \cstar) if $A$ is an inductive limit 
of \cstars\ of the form
 $$
 \bigoplus_{i=1}^ke_i(C_i\otimes 
 M_{n_i})e_i\oplus\bigoplus_{j=1}^rD_{m_j}\otimes M_{k_j},
 $$
where the first term is given as above. Some classes of ASH 
algebras are classified in terms of K theory (see, e.g., 
\cite{Ell,DG}); in particular a certain class of ASH algebras of 
real rank zero is classified by Dadarlat and Gong \cite{DG}.  The 
following follows from \ref{C5}, \ref{drop}, and \ref{A11} (3) and 
(5). 

\begin{theo} \label{C8}
Any ASH algebra (in the above sense) has Property \ref{A10}. 
\end{theo}

\section{Crossed products of AF algebras by $\Z$} 
 \setcounter{theo}{0} 
 \ncm{\cross}{{A\times_{\alpha}\Z}}      
Let $A$ be  an AF \cstar.  If $\alpha\in \Ap(A)=\AI(A)$  has the 
Rohlin property, then $\cross$ is a simple AT algebra and hence, 
as being an AH algebra,  satisfies Property \ref{A10}. More 
generally we have the following: 

\begin{theo}\label{F1}
Let $A$ be an AF \cstar\ and $\alpha\in\Aut(A)$. If the crossed 
product  $A\times_{\alpha}\Z$ is simple, then   
$A\times_{\alpha}\Z$ has Property \ref{A10}. 
\end{theo} 
\begin{pf} 
If   $A\times_{\alpha}\Z$ is isomorphic to the compact operators, 
then this follows from \ref{C8}. Hence we may assume that 
$A\times_{\alpha}\Z$ is not of type I.  (This assumption will be 
used only at the end of the proof.) 

Since $A$ is an AF \cstar, there is an increasing sequence $(A_n)$ 
of finite-dimensional $C^*$-subalgebras of $A$ such that 
$A=\ol{\cup_nA_n}$. By passing to a subsequence of $(A_n)$ we find 
a sequence $(u_n)$ in $\U(A)$ such that $u_{2n+1}\in A\cap 
A_{2n}'$ with $A_0=0$, $u_{2n}\in A\cap \Ad(u_{2n-1}u_{2n-2}\cdots 
u_1) \alpha(A_{2n-1})'$, $\|u_n-1\|<2^{-n}$, and
 $$
 \Ad\,u_1\alpha(A_1)\subset A_2\subset \Ad\,u_3u_2u_1\alpha(A_3)\subset
 A_4\subset \cdots.
 $$
By replacing $\alpha$ by $\Ad\,u\alpha$ with 
$u=\lim_nu_nu_{n-1}\cdots u_1$ and passing to a subsequence of 
$(A_n)$, we assume that
 $$
 \alpha^{\pm1}(A_n)\subset A_{n+1}
 $$
for all $n$. We denote by $1_n$ the unit of $A_n$. It is 
well-known that the crossed product $A\times_{\alpha}\Z$ remains 
the same for this inner perturbation of $\alpha$. 

Let $u$ denote the canonical unitary in $A\times_{\alpha}\Z$, 
i.e., $\cross$ is generated by $A$ and $u$ with relation that 
$ux=\alpha(x)u,\ x\in A$.

 Let $\F$ be a finite subset of $\cross$ and $\eps>0$. 
By taking a smaller $\eps>0$ if necessary, we can assume that $\F$ 
is $\{e_{s,ij}^{(k)}\ |\ i,j,s\}\cup \{1_ku1_k\}$, where 
$(e_{s,ij}^{(k)})$ is a family of matrix units for $A_k$ and 
$s\in\{ 1,2,\ldots,N_k\}$ corresponds to each direct summand of 
$A_k$. 

Assume that $n\in\N$ is so large that $2/n<\eps$. We take for $\G$ 
the set 
 $$\{1_ke_{s,ij}^{(k+2n)}u^me_{t,ab}^{(k+2n)}1_k\ |\ s,i,j,t,a,b; \ 
 m=0,\pm1,\pm2,\ldots,\pm 2n \}, 
 $$ 
where we assume that $1_k$ is the sum of some of 
$e_{t,aa}^{(k+2n)}$.  Let $B$ be a \cstar\ such that $B\supset 
\cross$ and let $\pi$ be an irreducible representation of $B$ with 
two unit vectors $\xi,\eta$ which define pure states $\om_1,\om_2$ 
of $B$ respectively. Let $\delta'>0$  and suppose that
 $$
 |\lan\pi(x)\xi,\xi\ran-\lan\pi(x)\eta,\eta\ran|<\delta',\ \ x\in\G.
 $$
We may assume that $\delta'$ is so small  that either, 
$\|1_k\xi\|<\eps$ and $\|1_k\eta\|<\eps$, or otherwise 
 $$
 |\lan u^me_{s,ij}^{(k+2n)}\xi_1,u^{\ell}e_{t,ab}^{(k+2n)}\xi_1\ran-
  \lan u^me_{s,ij}^{(k+2n)}\eta_1,u^{\ell}e_{t,ab}^{(k+2n)}\eta_1\ran|<\delta
  $$
for $m,\ell=0,1,2,\ldots,n$ and for some prescribed $\delta>0$, 
where we have omitted $\pi$, and $\xi_1=\|1_k\xi\|^{-1}1_k\xi$ and 
$\eta_1=\|1_k\eta\|^{-1}1_k\eta$.  

In the former case we should just take a continuous path $(v_t)$ 
in $\U((1-1_k)B(1-1_k))$ such that 
$v_1(1-1_k)\xi\approx(1-1_k)\eta$ and set $u_t=v_t+1_k$. Since 
$u_t$ commutes with $\F$ and $u_1\xi\approx\eta$, this completes 
the proof. 

In the latter case suppose that the linear space $\calL_{\xi_1}$ 
spanned  by $u^me_{s,ij}^{k+2n}1_k\xi_1$ with $m=0,1,\ldots,n$ and 
all possible $s,i,j$ is orthogonal to the space $\calL_{\eta_1}$ 
spanned by vectors of the same form with $\eta_1$ in place of 
$\xi_1$. If $\delta$ is sufficiently small, then we obtain vectors 
$\zeta(m,s,i,j)$ in $\calH_{\pi}\ominus L_{\xi_1}$ for 
$m=0,1,\ldots,n$ and for all $s,i,j$ with 
$e_{s,i,j}^{(k+2n)}1_k\neq0$ such that 
 $$
  \lan u^me_{s,ij}^{(k+2n)}\xi_1,u^{\ell}e_{t,ab}^{(k+2n)}\xi_1\ran=
  \lan\zeta(m,s,i,j),\zeta(\ell,t,a,b)\ran
  $$
and
 $$
 \|u^me_{s,ij}^{(k+2n)}1_k\eta-\zeta(m,s,i,j)\|<\eps'
 $$
for a small $\eps'$. Then we can define a projection $E$ on 
$\calL_{\xi_1}+\calL_{\zeta}$ such that
 \BE
 E( u^me_{s,ij}^{(k+2n)}1_k\xi+\zeta(m,s,i,j))&=&0,\\
 E(u^me_{s,ij}^{(k+2n)}1_k\xi-\zeta(m,s,i,j))
 &=&u^me_{s,ij}^{(k+2n)}1_k\xi-\zeta(m,s,i,j).   
 \EE
Then we find an $h=h^*\in B$ such that $\|h\|=1$, and 
$\pi(h)|\calL_{\xi}+\calL_{\zeta}=E$, where $\zeta$ denotes 
$\sum_{s,i}\zeta(0,s,i,i)$ and $\calL_{\zeta}$ is the space 
spanned by $\zeta(m,s,i,j)$'s. Define 
 $$
 \ol{h}=\frac{1}{n}\sum_{j=1}^{n-1}
     \sum_s\sum_iu^{-j}e_{s,i1}^{(k+n)}he_{s,1i}^{(k+n)}u^{j}.
 $$
Then  $\|[u,\ol{h}]\|<2/n<\eps$. Since 
$\sum_s\sum_ie^{(k+n)}_{s,i1}he^{(k+n)}_{s,1i}\in A\cap A_{k+n}'$ 
and $u^jA_ku^{-j}\subset A_{k+|j|}$, we have that $\ol{h}\in B\cap 
A_k'$. Also it follows that 
 \BE
 u^{-j}e_{s,i1}^{(k+n)}he_{s,1i}^{(k+n)}u^{j}(\xi_1+\zeta)&=&
 u^{-j}e_{s,i1}^{(k+n)}hu^j\alpha^{-j}( 
 e_{s,1i}^{(k+n)})(\xi_1+\zeta)  \\
 &\approx&u^{-j}e_{s,i1}^{(k+n)}h(\sum 
 c_{s,1i;t,ab}(u^{j}e^{(k+2n)}_{t,ab}\xi_1+\zeta(j,t,a,b))
 =0,
 \EE
where $\alpha^{-j}(e_{s,1i}^{(k+n)})1_k=\sum 
c_{s,1i;t,ab}e_{t,ab}^{(k+2n)}1_k$.  Here we have used the fact 
that $u^je_{t,ab}^{(k+2n)}\zeta\approx  
u^je_{t,ab}^{(k+2n)}\eta_1\approx \zeta(j,t,a,b)$. In the same way 
we obtain that 
 \BE
 u^{-j}e_{s,i1}^{(k+n)}he_{s,1i}^{(k+n)}u^{j}(\xi_1-\zeta)&=&
 u^{-j}e_{s,i1}^{(k+n)}hu^j\alpha^{-j}( 
 e_{s,1i}^{(k+n)})(\xi_1-\zeta)  \\
 &\approx&u^{-j}e_{s,i1}^{(k+n)}(\sum 
 c_{s,1i;t,ab}(u^{j}e^{(k+2n)}_{t,ab}\xi_1-\zeta(j,t,a,b))\\
 &\approx&u^{-j}e_{s,i1}^{(k+n)}u^j\alpha^{-j}(e_{s,1i}^{(k+n)})(\xi_1-\zeta)\\
 &=& \alpha^{-j}(e_{s,ii}^{(k+n)})(\xi_1-\zeta).
   \EE
Then, since $\alpha^{-j}(1_{k+n})\geq 1_k$, it follows that    
 \BE
 \ol{h}(\xi_1+\zeta)&\approx&0\\ 
 \ol{h}(\xi_1-\zeta)&\approx&\xi_1-\zeta.
 \EE
Hence 
 $$
 e^{i\pi\ol{h}}\xi_1=e^{i\pi\ol{h}}(\xi_1+\zeta)/2+e^{i\pi\ol{h}}(\xi_1-\zeta)/2 
 \approx (\xi_1+\zeta)/2 -(\xi_1-\zeta)/2
 \approx \zeta.
 $$
Thus we have that $e^{i\pi\ol{h}}1_k\xi\approx1_k\eta$.         
Note that $w_t=e^{i\pi t\ol{h}}$ in $\U(1_kB 1_k)$ almost commutes 
with $\F$. On the other hand we take a continuous path $(v_t)$ in 
$\U((1-1_k)B(1-1_k))$ such that $v_1(1-1_k)\xi\approx(1-1_k)\eta$. 
Taking the sum of $w_t$ and $v_t$ completes the proof in the case 
$\calL_{\xi} \perp \calL_{\eta}$.

If $\calL_{\xi} \not\perp \calL_{\eta}$, then we find a unit 
vector $\eta'$ such that $\calL_{\xi}\perp\calL_{\eta'},\ 
\calL_{\eta}\perp\calL_{\eta'}$, and
 $$
 |\lan x\eta,\eta\ran-\lan x\eta',\eta'\ran|<\delta',\ \ x\in\G
 $$
for an arbitrarily small $\delta'>0$. We apply the previous 
argument to the pairs $\xi,\eta'$ and $\eta',\eta$ to produce 
appropriate continuous paths $(u_t),(v_t)$ in $\U(B)$; in 
particular $u_1\xi=\eta'$ and $v_1\eta'=\eta$. Then the product 
$(v_tu_t)$ satisfies the required properties. 

To find such an $\eta'$ we use  the fact that $A\times_{\alpha}\Z$ 
is not of type I. Note that the set of vector states of 
$A\times_{\alpha}\Z$ in this representation 
$\rho=\pi|A\times_{\alpha}\Z$ is weak$^*$-dense in the state space 
of $A\times_{\alpha}\Z$. Hence there is a state $\varphi$ of 
$A\times_{\alpha}\Z$ such that $\pi_{\varphi}$ is disjoint from 
$\rho$ and 
 $$
 |\om_{\eta}\rho(x)-\varphi(x)|<\delta',\ \ x\in\G,
 $$                                                
where we denote by $\om_{\eta}$ the state on 
$\mathcal{B}(\calH_{\pi})$ defined by the vector $\eta$.  Then we 
can find a sequence $(\zeta_n)$ of unit vectors in $\calH_{\pi}$ 
such that $\om_{\zeta_n}\rho$ converges to $\varphi$ in the 
weak$^*$ topology. Hence $(\zeta_n)$ must converge to zero in the 
weak topology; we can take the required $\eta'$ near $\zeta_n$ for 
a sufficiently large $n$. 
\end{pf} 

From the proof of the above theorem we obtain:

\begin{theo} \label{F2}
Let $A$ be an AF \cstar\ and $\alpha\in\Aut(A)$. Then the crossed 
product $A\times_{\alpha}\Z$ has the transitivity: For any pair of 
pure states $\om_1$ and $\om_2$ of $A\times_{\alpha}\Z$ with 
$\ker\pi_{\om_1}=\ker\pi_{\om_2}$ there is an 
$\beta\in\AI(A\times_{\alpha}\Z)$ such that $\om_1=\om_2\beta$. 
\end{theo} 
\begin{pf}
If $\om_1\sim\om_2$ then this follows from Kadison's transitivity. 
Thus we may assume that $\om_1$ and $\om_2$ are not equivalent. 
Hence at least the quotient $A\times_{\alpha}\Z/\ker\pi_{\om_i}$ 
does not contain a non-zero type I ideal. We can prove this 
theorem just as Theorem \ref{A5} since $A\times_{\alpha}\Z$ has 
the following property (where $A$ denotes $A\times_{\alpha}\Z$). 
\end{pf}
 
\begin{property}  \label{F3}
For any finite subset $\F$ of $A$ and $\eps>0$ there exist a 
finite subset $\G$ of $A$ and $\delta>0$ satisfying: Let $B$ be a 
\cstar\ such that $B\supset A$.  For any pair of pure states 
$\om_1$ and   $\om_2$ of $B$ such that  $\om_1\sim\om_2$, 
$B/\ker\pi_{\om_1}$ does not contain a non-zero type I ideal, and 
 $$
 |\om_1(x)-\om_2(x)|<\delta,\ \ x\in \G,
 $$
there is a continuous path $(u_t)_{t\in[0,1]}$ in $\U(B)$ such 
that $u_0=1$, $\om_1=\om_2\Ad\,u_1$, and 
 $$
 \|\Ad\,u_t(x)-x\|<\eps,\ \ x\in\F.
 $$
\end{property}
 
\begin{lem}
Let $A$ be an AF \cstar\ and $\alpha\in\Aut(A)$. Then the crossed 
product $A\times_{\alpha}\Z$ has Property \ref{F3}. \end{lem} 
\begin{pf}
The proof of this fact immediately follows from the proof of 
Theorem \ref{F1}, since the assumption on $B/\ker\pi_{\om_1}$ 
substitutes the condition on $A\times_{\alpha}\Z$ there. \end{pf} 
 
\section{Purely infinite \cstars}
 \setcounter{theo}{0} 

Let $\lambda\in(0,1)$ and let $G_{\lambda}$ be the subgroup of 
$\R$ generated by $\lambda^n$ with $n\in\Z$. Equipping 
$G_{\lambda}$ with the order coming from $\R$, $G_{\lambda}$, 
being dense in $\R$, is a dimension group, i.e., there is a stable 
AF \cstar\ $A_{\lambda}$ whose dimension group is $G_{\lambda}$; 
$A_{\lambda}$ is unique up to isomorphism. Let $m_{\lambda}$ 
denote the automorphism of $G_{\lambda}$ defined by the 
multiplication of $\lambda$; there is an automorphism 
$\alpha_{\lambda}$ of $A_{\lambda}$ which induces $m_{\lambda}$ on 
$G_{\lambda}$; $\alpha_{\lambda}$ is unique up to cocycle 
conjugacy \cite{EK1}. If we denote by $\tau$ a (densely-defined 
lower-semicontinuous) trace on $A_{\lambda}$ (which is unique up 
to constant multiple), $\alpha_{\lambda}$ satisfies that 
$\tau\alpha_{\lambda}=\lambda\tau$. Then it follows from 
\cite{Ror} that the crossed product 
$A_{\lambda}\times_{\alpha_{\lambda}}\Z$ is a purely infinite 
simple \cstar. We can compute 
$K_*(A_{\lambda}\times_{\alpha_{\lambda}}\Z)$ by using the 
Pimsner-Voiculescu exact sequence; 
$K_0=G_{\lambda}/(1-\lambda)G_{\lambda}$ and $K_1=0$. When 
$\{f\in\Z[t]\ |\ f(\lambda)=0\}=p(t)\Z[t]$ with some $p(t)\in 
\Z[t]$, we have that $K_0=\Z/p(1)\Z$. By using the classification 
theory \cite{KP}, we know that  we get all the stable Cuntz 
algebras $\mathcal{O}_n\otimes \mathcal{K}$ in this way. Recall 
that $K_0(\mathcal{O}_n)=\Z/(n-1)\Z$ if $n<\infty$ and 
$K_0(\mathcal{O}_{\infty})=\Z$. Hence by using \ref{A11} and 
\ref{C5} we have: 

\begin{lem}
All the \cstars\ of the form
 $$
 C(\T)\otimes M_k\otimes\mathcal{O}_n
 $$
have Property \ref{A10}, where $n=2,3,\ldots,\infty$ with 
$n=\infty$ inclusive and $k=1,2,\ldots,n-1$ if $n<\infty$ and 
$k=1,2,\ldots$ otherwise. Moreover all the \cstars\ obtained as 
inductive limits of finite direct sums of \cstars\ of the above 
form have Property \ref{A10}. 
\end{lem} 

Note that $K_0(C(\T)\otimes M_k\otimes\mathcal{O}_n) 
=K_0(\mathcal{O}_n)$ and $K_1(C(\T)\otimes 
M_k\otimes\mathcal{O}_n) =K_0(\mathcal{O}_n)$. Tensoring with 
$M_k$ exhausts all possible pairs $(K_0,[1])$ for cyclic groups 
$K_0$. In \cite{BEEK4} a class of inductive limit \cstars\ is 
considered using the \cstars\ in the above lemma as building 
blocks. It is not difficult to show that the purely infinite 
simple separable unital \cstars\ obtained this way exhaust all 
possible $(G_0,g,G_1)$ for $(K_0,[1],K_1)$, where $G_0$ and $G_1$ 
are arbitrary countable abelian groups and $g\in G_0$. (Given any 
pair of  countable abelian groups $G_0,G_1$ and $g\in G_0$ we find 
an inductive system $(G_{in},\phi_{in})$ whose limit is $G_i$, 
where all $G_{in}$ is a finite direct sum of cyclic groups; for 
$i=0$ we specify $g_{0n}\in G_{0n}$ with 
$\phi_{0n}(g_{0n})=g_{0,n+1}$ so that $(g_{0n})$ represents $g$. 
Let $G_n=G_{0n}\oplus G_{1n}$ and $g_n=g_{0n}\oplus0$ and extend 
$\phi_{in}:G_{in}\ra G_{i,n+1}$ to a map $\phi_{in}:G_n\ra 
G_{n+1}$ by adding zero maps. Then it follows that $G_i$ is the 
inductive limit of $(G_n,\phi_{in})$ for $i=0,1$. For each  $G_n$ 
we take a direct sum $A_n$ of \cstars\ of the form  $C(\T)\otimes 
M_k\otimes\mathcal{O}_n$ with $K_*(A_n)=G_n$ and $[1]=g_n$. Then 
we find a unital homomorphism $\varphi_n:A_n\ra A_{n+1}$ such that 
$\varphi_n$ induces $\phi_{in}$ for $i=0,1$ and $\varphi_n$ is a 
non-zero map from each direct summand of $A_n$ into each direct 
summand of $A_{n+1}$; moreover $\varphi_n$ maps the canonical 
unitary $z\otimes 1\otimes 1$ of each direct summand  
$C(\T)\otimes M_k\otimes\mathcal{O}_n$ of $A_n$ to a non-zero 
unitary $u$ for each direct summand of $A_{n+1}$ such that the 
spectrum of $u$ evaluated at each $t\in\T$ is $\T$ (a considerably 
weaker condition will suffice; see \cite{BEEK4} for details); the 
latter is imposed to insure that the inductive limit $A$ of 
$(A_n,\varphi_n)$ is simple. Then $A$ satisfies that $K_i(A)=G_i$ 
and $[1]=g$.)

Let $\calC_{\infty}$ denote the class of  purely infinite 
separable simple \cstars\ satisfying the Universal Coefficient 
Theorem, i.e., $\calC_{\infty}$ is the class of \cstars\ 
classified in terms of K theory by Kirchberg and Phillips 
\cite{KP,Phi}.  Since all of $\calC_{\infty}$ can be obtained as 
inductive limits as discussed above, we have:

\begin{theo}
Any \cstar\ $A$ in  $\calC_{\infty}$ has Property \ref{A10}. In 
particular  $\AI(A)$ acts transitively on $P(A)$. 
\end{theo}

\section{The group \cstars\ of discrete amenable groups}
 \setcounter{theo}{0}
\begin{theo}
If $G$ is a countable discrete amenable group, then the group 
\cstar\ $C^*(G)$ satisfies the transitivity: for any pair of pure 
states $\om_1$ and $\om_2$ with 
$\ker{\pi_{\om_1}}=\ker{\pi_{\om_2}}$ there is an 
$\alpha\in\AI(C^*(G))$ such that $\om_1=\om_2\alpha$. 
\end{theo} 
\begin{pf}
If $\om_1\sim\om_2$ then this follows from Kadison's transitivity.
Thus we may assume that $\om_1$ and $\om_2$ are not equivalent. 
Hence at least the quotient $C^*(G)/I$ with $I=\ker\pi_{\om_i}$ 
does not contain a non-zero type I ideal. We can prove this 
theorem just as Theorem \ref{A5} (or \ref{F2}) once we have shown 
the following: 
\end{pf}
 
\begin{lem} 
If $G$ is a  discrete amenable group, $C^*(G)$ has Property 
\ref{F3}. 
\end{lem}
\begin{pf}
Since G is a discrete amenable group, for any finite subset $\F$ 
of $G$ and $\eps>0$ there is a finite subset $\G$ of $G$ such that 
 $$\frac{|\G\vartriangle\G g|}{|\G|}<\eps,\ \ g\in \F,
 $$
where $\vartriangle$ denotes the difference of sets and 
$|\,\cdot\,|$ denotes the number of elements. We may suppose that 
$\G\ni 1$. 

 Suppose that $C^*(G)$ is a unital $C^*$-subalgebra of 
$B$ and let $\pi$ be an irreducible representation of $B$. Let 
$\delta>0$ and let $\xi,\eta$ be unit vectors in $\calH_{\pi}$ 
such that  
 $$
 |\lan\pi(g)\xi, \pi(h)\xi\ran
 -\lan \pi(g)\eta, \pi(h)\eta 
 \ran|<\delta  
 $$
for $g,h\in \G$. 

First suppose that the linear subspace $\calL_{\xi}$  spanned by  
$\pi(g)\xi,\ g\in\G$ and the linear subspace $\calL_{\eta}$ 
spanned by $\pi(g)\eta, \ g\in\G$ are mutually orthogonal. Then 
for a sufficiently small $\delta>0$ there is a family $(\zeta(g))$ 
of vectors in $\calH\ominus\calL_{\xi}$, which is 
infinite-dimensional, such that             
 $$
 \|\pi(g)\eta-\zeta (g)\|<\eps', \ \ g\in\G,
  $$
for small $\eps'>0$ and 
 $$
 \lan\pi(g)\xi, \pi(h)\xi\ran
 =\lan\zeta(g),\zeta(h)\ran. 
 $$ 
Since $\pi(g)\zeta\approx\pi(g)\eta\approx\zeta(g)$ for $g\in\G$ 
with $\zeta=\zeta(1)$, we have that 
$\|\pi(g)\zeta-\zeta(g)\|<2\eps'$.  Then there is a projection $E$ 
on $\calL_{\xi}+\calL_{\zeta}$ such that 
 \BE  E( \pi(g)\xi+ 
 \zeta(g))&=&0,\\
     E( \pi(g)\xi- \zeta(g))
     &=&\pi(g)\xi- \zeta(g).
 \EE
where $\calL_{\zeta}$ is the subspace spanned by $\zeta(g),\ 
g\in\G$. Then we find an $h=h^*\in B$ such that $\|h\|=1$ and 
 $$
 \pi(h)|(\calL_{\xi}+\calL_{\zeta})=E.
 $$
Set 
 $$\ol{h}=\frac{1}{|\G|}\sum_{g\in\G}g^{-1}hg,
 $$  
which is a self-adjoint element in $B$ with norm at most one and 
satisfies that $\|g^{-1}\ol{h}g-\ol{h}\|\leq \eps,\ g\in\F$. Note 
that for $g\in\G$, 
 $$
 \|\pi(g^{-1}hg)(\xi+\zeta)\|<2\eps'
 $$
and
 $$\|\pi(g^{-1}hg)(\xi-\zeta)-(\xi-\zeta)\|<2\eps'.
 $$
Then we obtain that $\pi(e^{i\pi \ol{h}})\xi\approx\zeta$ because
 \BE
 &&\|\pi(e^{i\pi \ol{h}})\xi-\zeta\| \\
 &&\leq 
 \|\sum_{m=1}^{\infty}\frac{(i\pi)^m}{m!}\pi(\ol{h}^m)\frac{1}{2}(\xi+\zeta)\|
 +\|\sum_{m=1}^{\infty}\frac{(i\pi)^m}{m!}(\pi(\ol{h}^m)-1) \frac{1}{2}(\xi-\zeta)\| 
 \\
 &&\leq \eps'(e^{\pi}-1+\pi e^{\pi}).
 \EE
Define $u_t=e^{it\pi \ol{h}}$, which is a continuous path in 
$\U(B)$ satisfying $\|[u_t,g]\|<\pi\eps,\ g\in\F$ and 
$\pi(u_1)\xi\approx\zeta\approx \eta$. Hence this completes the 
proof in the case that $\calL_{\xi}\perp\calL_{\eta}$.  
          
If  $\calL_{\xi}\not\perp\calL_{\eta}$, then we first find a unit 
vector $\eta'$ such that $\calL_{\xi}\perp\calL_{\eta'}$, 
$\calL_{\eta}\perp\calL_{\eta'}$, and 
 $$
  |\lan\pi(g)\eta, \pi(h)\eta\ran
 -\lan \pi(g)\eta', \pi(h)\eta' 
 \ran|<\delta',\ \ g,h\in\G
 $$         
for a very small $\delta'>0$. We apply the previous argument to 
the pairs $\xi,\eta'$ and $\eta',\eta$ to get the conclusion. 

To find such an $\eta'$ we may argue as in the proof of \ref{A3} 
(see also the final part of the proof of \ref{F1}). Let $e\in B$ 
be such that $0\leq e\leq 1$ and 
 $$
 \|e h^{-1}g e-\om_2(\pi(h^{-1}g))e^2\|<\delta',\ \ g,h\in\G,
 $$                                                          
where $\om_2$ is the vector state defined by $\eta$. If there is 
no such $\eta'$, then the spectral projection of $\pi(e)$ 
corresponding to $[1-\delta',1]$ must be finite-dimensional; i.e., 
$\pi(B)$ contains the compact operators on $\calH_{\pi}$, which is 
excluded by the assumption.                     
\end{pf}

\section{Strong transitivity}
 \setcounter{theo}{0}
We first introduce the following property for a \cstar\ $A$ which 
is stronger than \ref{A10}:

\begin{property}\label{S1} 
For any finite subset $\F$ of $A$ and $\eps>0$ there exist a 
finite subset $\G$ of $A$ and $\delta>0$ satisfying: If $B$ is a 
\cstar\ containing $A$ as a $C^*$-subalgebra and 
$\omega_1,\om_2,\ldots,\om_n,\varphi_1,\varphi_2,\ldots,\varphi_n$ 
are pure states of $B$ such that 
 \begin{enumerate}
 \item $(\om_i)$ are mutually disjoint,
 \item $\om_i\sim\varphi_i$ for $i=1,2,\ldots,n$,
 \item $ |\omega_i(x)-\varphi_i(x)|<\delta,\ \ x\in\G,\ i=1,2,\ldots,n$,
 \end{enumerate}
then  there is a continuous path $(u_t)_{t\in[0,1]}$ in $\U(B)$ 
such that $u_0=1$, $\omega_i=\varphi_i\Ad\,u_1,\ i=1,2,\ldots,n$, 
and 
 $$ 
 \|\Ad\,u_t(x)-x\|<\eps,\ \ x\in\F,\ t\in[0,1].
 $$
\end{property} 

We can actually show the above property for all the \cstars\ for 
which we have shown Property \ref{A10}. This is because we can use 
the same argument for each pair $\om_i,\varphi_i$ up to the point 
where we invoke Kadison's transitivity; instead of working in one 
irreducible representation of $B$ we are now working in  
$\bigoplus_{i=1}^n\pi_{\om_i}$, which is a direct sum of mutually 
disjoint irreducible representations. We have to find an element 
$h$ in the \cstar\ $B$ with the prescribed property in this larger 
space. But there is a form of Kadison's transitivity in this 
generality (cf. 1.21.16 of \cite{Sak}); so we are done. 

The following follows just like \ref{A3}.

\begin{lem}
If $\om_1,\om_2,\ldots,\om_n,\varphi_1,\varphi_2,\ldots,\varphi_n$ 
are pure states of $A$ such that 
 \begin{enumerate}
 \item $(\varphi_i)$ are mutually disjoint,
 \item $\ker\pi_{\om_i}=\ker\pi_{\varphi_i}$ for $i=1,2,\ldots,n$,
 \end{enumerate}
then for any finite subset $\F$ of $A$ and $\eps>0$ there is a 
$u\in\U(A)$ such that 
 $$|\om_i(x)-\varphi_i\Ad\,u(x)|<\eps,\ \ x\in\F,\ i=1,2,\ldots,n.
 $$
\end{lem}   
  
\begin{theo}\label{S2}
Suppose that $A$ is a separable  \cstar\ satisfying Property 
\ref{S1}. Let $(\om_i)_{1\leq i\leq n}$ and $(\varphi_i)_{1\leq 
i\leq n}$ be finite sequences of pure states of $A$ such that 
$(\om_i)$ (resp. $(\varphi_i)$) are mutually disjoint and 
$\ker\pi_{\om_i}=\ker\pi_{\varphi_i}$ for all $i$. Then there is 
an $\alpha\in\AI(A)$ such that $\om_i=\varphi_i\alpha$ for all 
$i=1,2,\ldots,n$. In particular if $(\om_i)_{1\leq i\leq n+1}$ are 
pure states such that  they are mutually disjoint maybe except for 
the pair $\om_1,\om_{n+1}$ and all $\ker\pi_{\om_i}$'s are equal, 
then there is an $\alpha\in\AI(A)$ such that 
$\om_i\alpha=\om_{i+1},\ i=1,2,\ldots,n$.   
\end{theo}
   
The proof goes just like the proof of \ref{A5} does. We present 
what we use at each induction step in a form of lemma:

\begin{lem}\label{S3}
Suppose that $A$  satisfies \ref{S1}.  For any finite subset $\F$ 
of $A$ and $\eps>0$ there exist a finite subset $\G$ of $A$ and 
$\delta>0$ satisfying: If $\omega_i,\varphi_{i},\ i=1,2,\ldots,n,$ 
are pure states of $A$ such that 
 \begin{enumerate}
 \item $(\varphi_i)$ are mutually disjoint,
 \item $\ker\pi_{\om_i}=\ker\pi_{\varphi_i}$ for all $i$,
 \item $ |\omega_i(x)-\varphi_i(x)|<\delta,\ \ x\in\G,\ i=1,2,\ldots,n$,
 \end{enumerate}
then  for any finite subset $\F'$ of $A$ and $\eps'>0$ there is a 
continuous path $(u_t)_{t\in[0,1]}$ in $\U(A)$ such that $u_0=1$, 
 \BE
   |\omega_i(x)-\varphi_i\Ad\,u_1(x)|&<&\eps',\ \ x\in\F', \  i=1,2,\ldots,n,\\ 
  \|\Ad\,u_t(x)-x\|&<&\eps,\ \ x\in\F,\ t\in[0,1].
 \EE
\end{lem} 

Finally we present another version of transitivity:

\begin{theo}\label{S4}
Suppose that $A$ is a separable \cstar\ with Property \ref{S1}. 
Let $(\pi_n)$ and $(\rho_n)$ be sequences of irreducible 
representations of $A$ such that $(\pi_n)$ (resp. $(\rho_n)$) are 
mutually disjoint and $\ker\pi_n=\ker\rho_n$ for all $n$. Then 
there is an $\alpha\in\AI(A)$ such that $\pi_n=\rho_n\alpha$ for 
all $n\in\N$. In particular if $(\pi_n)_{n\in {\bf Z}}$ are 
irreducible representations of $A$ such that they are mutually 
disjoint and all $\ker\pi_n$'s are equal, then there is an 
$\alpha\in\AI(A)$ such that $\pi_n\alpha=\pi_{n+1}$ for all $n$. 
\end{theo} 
\begin{pf}   
The proof is similar to the one of \ref{S2}. We will construct 
pure states $\om_n$ (resp. $\varphi_n$) associated with $\pi_n$ 
(resp. $\rho_n$) inductively such that $\om_n=\varphi_n\alpha$ 
holds for all $n$. To introduce a new pair of pure states at each 
induction step we will use the following easy lemma. 
\end{pf}

\begin{lem}
Let $\pi$ and $\rho$ be mutually disjoint irreducible 
representations of $A$ with $\ker\pi=\ker\rho$ and $u\in\U(A)$. 
Then for any finite subset $\G$ of $A$ and $\delta>0$ there are a 
pure state $\om$ associated with $\pi$ and a pure state $\varphi$ 
associated with $\rho$ such that 
$|\om(x)-\varphi\Ad\,u(x)|<\delta,\ x\in\G$. 
\end{lem} 
   
\begin{rem}
If $A$ is a non type I separable simple \cstar\ with Property 
\ref{S1}, it follows from the above theorem that there is an 
$\alpha\in\AI(A)$ such that all non-zero powers $\alpha^n$ are 
outer. Since for any injective map $\om$ of $\Z$ into the set of 
equivalence classes of irreducible representations of $A$ there is 
an $\alpha\in\AI(A)$ such that $\om(n)\alpha=\om(n+1),\ n\in\Z$, 
it follows that the quotient $\AI(A)/\Inn(A)$ is uncountable. 
\end{rem}

\section{Remarks on Property \ref{A10}}
 \setcounter{theo}{0}  
 
Looking at the proofs of \ref{C4}--\ref{drop} and \ref{F1} etc., 
we come to know that the path $(u_t)$ in Property \ref{A10} may be 
chosen so that its length is dominated by a universal constant, 
which is only slightly bigger than $\pi$. (This follows by 
modifying the proofs given there; it is $2\pi$ instead of $\pi$ 
that follows immediately.) Taking  this fact  into consideration, 
we first introduce the following stronger condition: 
 
 \ncm{\length}{{\rm length}}
\begin{property}\label{R1}(for a \cstar\ $A$ and a constant $C\geq\pi$) 
For any finite subset $\F$ of $A$ and $\eps>0$ there exist a 
finite subset $\G$ of $A$ and $\delta>0$ satisfying: If $B$ is a 
\cstar\ containing $A$ as a $C^*$-subalgebra and  $\omega_1$ and 
$\omega_2$ are pure states of $B$ such that $\omega_1\sim 
\omega_2$ and 
 $$ |\omega_1(x)-\omega_2(x)|<\delta,\ \ x\in\G,
 $$
then for any $\eps'>0$ there is a rectifiable path 
$(u_t)_{t\in[0,1]}$ in $\U(B)$ such that $u_0=1$, 
$\omega_1=\omega_2\Ad\,u_1$, $\length((u_t))<C+\eps'$, and 
 $$ \|\Ad\,u_t(x)-x\|<\eps,\ \ x\in\F,\ t\in [0,1].
 $$
\end{property}

As asserted above, in the cases we handled in the previous 
sections, the \cstar\ $A$ has this property for $C=\pi$. If 
$A=\C$, then $C=\pi/2$ suffices for the above property to hold. If 
the property holds for all $A=\C^n$ with $n\in\N$ and for a 
constant $C$, then one can  check that $C\geq\pi$. To illustrate 
these points we give two propositions: 

\ncm{\Spec}{{\rm Spec}}
\begin{prop}
Let $\calH$ be a Hilbert space and $e$ be a projection in $\calH$. 
If $\xi$ and $\eta$ are unit vectors in $\calH$ such that $\lan 
e\xi,\xi\ran =\lan e\eta,\eta\ran$, there is a rectifiable path 
$(u_t)_{t\in [0,1]}$ in the group $\U(\calH)$ of unitaries on 
$\calH$ such that $u_0=1$, $u_1\xi=\lambda\eta$ for some 
$\lambda\in\T$, $[u_t,e]=0$, and $\length((u_t))\leq\pi/2$. 
\end{prop}
\begin{pf}
There is a $\lambda\in\T$ such that 
 $$
 \Re\lan e\xi,\lambda\eta\ran\geq0,\ \ 
 \Re\lan (1-e)\xi,\lambda\eta\ran\geq0.
 $$
Hence we only have to apply the following proposition to the pairs 
$e\xi,\lambda e\eta\in e\calH$ and $(1-e)\xi,\lambda(1-e)\eta\in 
(1-e)\calH$ separately to reach the conclusion. 
\end{pf}

\begin{prop}
Let $\calH$ be a  Hilbert space and $\xi,\eta$ unit vectors in 
$\calH$. If $\theta\in[0,\pi]$ is such that 
$\Re\lan\xi,\eta\ran=\cos\theta$, then there is a rectifiable path 
$(u_t)$ in $\calU(\calH)$ such that $u_0=1$, $u_1\xi=\eta$, and 
$\length((u_t))=\theta$. Moreover for any rectifiable path $(v_t)$ 
in $\calU(\calH)$ with $v_0=1$ and $v_1\xi=\eta$, it follows that 
$\length((v_t))\geq\theta$. 
\end{prop}
\begin{pf}
If $\eta=e^{\pm i\theta}\xi$, then we set $u_t=e^{\pm it\theta}$. 

Suppose that $\xi$ and $\eta$ are linearly independent. For 
$t\in[0,1]$ let 
 $$
 \xi(t)=\cos\theta t\cdot \xi+\sin\theta t(\sin\theta)^{-1}
                                       (\eta-\cos\theta\cdot\xi).
 $$
Then $\xi(0)=\xi$ and $\xi(1)=\eta$. Since 
$\Re\lan\xi,\eta-\cos\theta\cdot\xi\ran=0$ and 
$\|\eta-\cos\theta\cdot\xi\| =\sin\theta$, it follows that 
$\|\xi(t)\|=1$ and $\|\xi'(t)\|=\theta$. Thus $(\xi(t))$ is a 
$C^1$ path in the unit vectors of $\calH$ and its length is 
$\theta$. 

We define $\alpha(t)\in\R$ by 
 $$
 i\alpha(t)=\lan\xi'(t),\xi(t)\ran=i\theta(\sin\theta)^{-1}
           \Im\lan\eta,\xi\ran
 $$
and $\beta(t)\in\R$ by 
 $$
 \beta(t)=\|\xi'(t)-i\alpha(t)\xi(t)\|=(\theta-\alpha(t))^{1/2}.
 $$
Since both $\alpha(t)$ and $\beta(t)$ are constants, we write 
$\alpha(t)=\alpha$ and $\beta(t)=\beta$. Let 
$\zeta(t)=\beta^{-1}(\xi'(t)-i\alpha\xi(t))$ and define a 
self-adjoint operator $h(t)$ by 
 $$
 h(t)=\alpha\xi(t)\otimes\xi(t)-i\beta\zeta(t)\otimes\xi(t)
       +i\beta\xi(t)\otimes\zeta(t) -\alpha\zeta(t) \otimes\zeta(t).
 $$
Then $(h(t))$ is a continuous path in the set of self-adjoint 
operators in $\calH$ such that $ih(t)\xi(t)=\xi'(t)$ and 
 $$
 \|h(t)\|=(\alpha^2+\beta^2)^{1/2}=\theta=\|\xi'(t)\|.
 $$
We define a $C^1$ path $(u(t))_{t\in[0,1]}$ in $\calU(\calH)$ by 
 $$
 \frac{d}{dt}u(t)=ih(t)u(t),\ \ u(0)=1.
 $$
In fact, since $\frac{d}{dt}u(t)^*u(t)=0=\frac{d}{dt}u(t)u(t)^*$, 
we have that $u(t)^*u(t)=1=u(t)u(t)^*$. We also have that 
$\length((u(t))=\theta$. Since $u(0)\xi=\xi=\xi(0)$ and 
 \BE
 \frac{d}{dt}\|u(t)\xi-\xi(t)\|^2&=&\frac{d}{dt}(2-2\Re\lan 
 u(t)\xi,\xi(t)\ran)\\
 &=&-2\Re\lan ih(t)u(t)\xi,\xi(t)\ran-2\Re\lan u(t)\xi,\xi'(t)\ran\\
 &=&0,
 \EE
we have that $u(t)\xi=\xi(t)$; in particular $u(1)\xi=\eta$. Thus 
$(u(t))$ is the desired path in $\U(\calH)$. 

Let $(v_t)$ be a rectifiable path in $\U(\calH)$ such that $v_0=1$ 
and $v_1\xi=\eta$. Since 
 $$
 2\cos\theta=2\Re\lan\xi,\eta\ran=\lan(v_1^*+v_1)\xi,\xi\ran,
 $$
it follows that the spectrum $\Spec(v_1)$ has an $e^{i\varphi}$ 
such that $\cos\varphi\leq\cos\theta=\Re\lan\xi,\eta\ran$. For any 
$\eps>0$ we find  a sequence $(t_0,t_1,\ldots,t_n)$ in $[0,1]$ 
such that $t_0=0<t_1<t_2<\cdots<t_n=1$ and 
$\|v_{t_i}-v_{t_{i-1}}\|<\eps$. Then, by the following lemma, we 
find a $\lambda_i\in\Spec(v_{t_i})$ such that 
$\lambda_n=e^{i\varphi}$ and 
$|\lambda_i-\lambda_{i-1}|\leq\|v_{t_i}-v_{t_{i-1}}\|$. Hence 
 $$
 \sum_{i=1}^n|\lambda_i-\lambda_{i-1}|\leq\sum_{i=0}^n\|v_{t_i}-v_{t_{i-1}}\|
  \leq \length((v_t)).
 $$
Since $\lambda_0=1,\, \lambda_n=e^{i\varphi}$ and $\eps>0$ is 
arbitrary, we obtain that $\varphi\leq\length((v_t))$. Since 
$\theta\leq\varphi$, this completes the proof. 
\end{pf} 

\begin{lem}
Let $A$ be a \cstar\ and let $u,v\in\U(A)$. If 
$\lambda\in\Spec(u)$, then there is a $\mu\in\Spec(v)$ such that 
$|\lambda-\mu|\leq \|u-v\|$. 
\end{lem}
\begin{pf}
Let $\delta=\|u-v\|$ and $w=u^*v$. Since $\|1-w\|=\delta$, it 
follows that $\Spec(w)\subset \{e^{it}\ |\ 
|t|\leq\theta\}=\T_{\theta}$, where $\theta=2\sin^{-1}\delta/2$. 
Let $\om$ be a state of $A$ such that $\om(u)=\lambda\in\Spec(u)$. 
Then $\om(w)=\om(u^*v)=\ol{\lambda}\om(v)$. Hence 
$\lambda\om(w)=\om(v)$ belongs to the convex closure of 
$\lambda\T_{\theta}$, which implies that $\Spec(v)\cap 
\lambda\T_{\theta}\not=\emptyset$. 
\end{pf}

We seem to need the above stronger property to prove: 

\begin{prop}\label{R2}
If a \cstar\ $A$ has Property \ref{R1} for a constant $C\geq \pi$ 
and $A_1$ is a hereditary $C^*$-subalgebra of $A$, then $A_1$ has 
Property \ref{R1} for the same constant $C$. \end{prop} 

To show this we first present the non-unital version of 
\ref{A105}: 

\begin{lem}\label{R3}
When $A$ is non-unital, Property \ref{R1} is equivalent to the one 
obtained by restricting the ambient \cstar\ $B$ to a \cstar\ 
having an approximate identity for $A$ as an approximate identity 
for $B$ itself. \end{lem} 
\begin{pf}
Technically the proof will be similar to the proofs of \ref{mult} 
and \ref{drop}. 

We assume the weaker property for $A$: For $(\F,\eps)$ there is a 
$(\G,\delta)$ satisfying: If $B\supset A$ and $B=\ol{ABA}$, and 
$\om_1$ and $\om_2$ are pure states of $B$ such that 
$\om_1\sim\om_2$ and $|\om_1(x)-\om_2(x)|<\delta,\ x\in\G$, then 
for any $\eps'>0$ there is a rectifiable path $(u_t)_{t\in [0,1]}$ 
in $\U(B)$ such that $u_0=1$, $\om_1=\om_2\Ad\,u_1$, 
$\length((u_t))<C+\eps'$, and $\|\Ad\,u_t(x)-x\|<\eps,\ 
x\in\F,\,t\in [0,1]$.

Let $\F$ be a finite subset of $A$ and $\eps>0$. We may assume 
that $\|x\|\leq1$ for $x\in\F$, $\eps>0$ is sufficiently small, 
and that there is an $e\in A$ such that $0\leq e\leq 1$ and 
$exe=x,\ x\in\F$. Let $B$ be a \cstar\ with $B\supset A$. 
                                                 
Let $\eps'>0$ be such that $2\sqrt{\eps'}(C+8\sqrt{\eps'})<\eps$  
and let $M,N\in\N$ be so large that $M\eps'{^2}>4$ and 
$N\eps'>2C$. There is a sequence $(e_0,e_1,\ldots,e_{M(N+1)})$ in 
$A$ such that $e_0xe_0=x,\ x\in\F$, $0\leq e_i\leq 1$, and 
$e_ie_{i-1}=e_{i-1}$. Let $d>1$ be a very large constant and let 
$\F_1=\F\cup \{de_i\ | \ 0\leq i\leq M(N+1)\}$. For $(\F_1,\eps)$ 
we choose a $(\G,\delta)$ as in the weaker version of Property 
\ref{R1}. We may assume that $\|x\|\leq 1,\ x\in\G$ and 
$\delta<1$. 

Let $\pi$ be an irreducible representation of $B$ and $\xi,\eta$ 
unit vectors in $\calH_{\pi}$. Suppose that 
 $$
 |\lan\pi(x)\xi,\xi\ran-\lan\pi(x)\eta,\eta\ran|<\delta\eps'{^2}/2,\ \ x\in\G.
 $$
There exists a $j$ between $1$ and $(M-1)(N+1)$ inclusive such 
that 
 $$
 \lan\pi(e_{j+N}-e_{j-1})\xi,\xi\ran+\lan\pi(e_{j+N}-e_{j-1})\eta,\eta\ran
 <\eps'{^2}/2.
 $$
(Otherwise we would be led to a contradiction, $M\eps'{^2}/2\leq 
2$.) If $\lan\pi(e_{j-1})\xi,\xi\ran>\eps'{^2}$, then we have that  
for $\xi_1=\|\pi(e_{j-1}^{1/2})\xi\|^{-1}\pi(e_{j-1}^{1/2})\xi$ 
and 
$\eta_1=\|\pi(e_{j-1}^{1/2})\eta\|^{-1}\pi(e_{j-1}^{1/2})\eta$, 
 $$
 |\lan\pi(x)\xi_1,\xi_1\ran-\lan\pi(x)\eta_1,\eta_1\ran|<\delta,\ \ x\in\G.
 $$ 
Here we have used that 
 $$
 |1-\frac{\lan\pi(e_{j-1})\eta,\eta\ran} 
 {\lan\pi(e_{j-1})\xi,\xi\ran}|<\delta/2.
 $$
Then, since $\xi_1,\eta_1\in \pi(A)\calH_{\pi}$, there is a 
rectifiable path $(u_t)$ in $\U(\ol{ABA})$ such that $u_0=1$, 
$\pi(u_1)\xi_1=\eta_1$, $\length((u_t))<C+\eps'$, and 
$\|\Ad\,u_t(x)-x\|<\eps,\ x\in\F_1$. Moreover we may assume that 
the length of $(u_t)_{t\in [0,s]}$ is proportional to $s$ for any
$s\in[0,1]$. We set $p_k=e_{j+k}-e_{j+k-1}$ and define
 $$
 z_t=u_te_j+\sum_{k=1}^{N-1}u_{(1-k/N)t}p_k+1-e_{j+N-1}.
 $$                                                     
Let $B_0$ be the closure of $e_{j+N-1}Be_{j+N-1}$. Then $(z_t)$ is 
a path in $B_0+1$, from which we shall construct a path $(v_t)$ in 
$\U(B_0)$ with appropriate properties. 
 
Note that $\|[u_t,e_i]\|<\eps/d$, where $d\gg1$ is chosen 
independently of $N$; i.e., $\|[u_t,e_i]\|\approx0$. Since 
$z_tx=u_tx$ and $xz_t\approx xu_t$ for $x\in\F$, we have that 
$\|[z_t,x]\|\approx \|[u_t,x]\|<\eps$ for $x\in\F$. Thus by making 
$d$ sufficiently large, we may assume that $\|[z_t,x]\|<\eps,\ 
x\in\F$.

Let $p_0=e_j$ and $p_N=1-e_{j+N-1}$. Then $\sum_{k=0}^Np_k=1$ and 
$p_{k}p_{\ell}=0$ if $|k-\ell|>1$. Let $s(p_k)$ denote the support 
projection of $p_k$ in $B^{**}$. Then we have that
 \BE
 &&\|(z_t-u_{(1-k/N)t})s(p_k)\|^2\\
 &&=\|(u_{(1-(k-1)/N)t}p_{k-1}+u_{(1-k/N)t}p_k+u_{(1-(k+1)/N)t}p_{k+1}-u_{(1-k/N)t})s(p_k)\|^2\\
 &&=\|\{(u_{(1-(k-1)/N)t}-u_{(1-k/N)t})p_{k-1}
 +(u_{(1-(k+1)/N)t}-u_{(1-k/N)t})p_{k+1}\}s(p_k)\|^2\\  
 &&\approx 
 \|s(p_k)\{p_{k-1}(u_{(1-(k-1)/N)t}-u_{(1-k/N)t})^*(u_{(1-(k-1)/N)t}-u_{(1-k/N)t})p_{k-1}\\
 &&\ 
 +p_{k+1}(u_{(1-(k+1)/N)t}-u_{(1-k/N)t})^*(u_{(1-(k+1)/N)t}-u_{(1-k/N)t})p_{k+1}\}s(p_k)\|\\
 &&\leq (\frac{C+\eps'}{N})^2,
 \EE                          
where $p_{-1}=0=p_{N+1}$. Since $C/N<\eps'/2$, we may assume that 
 $$
 \|(z_t-u_{(1-k/N)t})s(p_k)\|<\eps'.
 $$
Since 
 \BE
 \|z_t^*z_t-1\|&\approx& \|\sum_{k=0}^Np_k^{1/2}(z_t^*z_t-1)p_k^{1/2}\|\\
 &=& \|\sum_{k=0}^Np_k^{1/2}\{(z_t^*-u_{(1-k/N)t}^*)z_t
 +u_{(1-k/N)t}^*(z_t-u_{(1-k/N)t})\}p_k^{1/2}\|\\
 &<&\eps'(2+\eps'),
 \EE
we may assume that $\|z_t^*z_t-1\|<3\eps'$. For $0\leq s<t\leq1$, 
let
 $$y(s,t)=\sum_{k=0}^Np_k^{1/2}(u_{(1-k/N)t}-u_{(1-k/N)s})^*(u_{(1-k/N)t}-u_{(1-k/N)s})p_k^{1/2}.
 $$
Then we have that
 $$
 \|y(s,t)\|\leq (C+\eps')^2(t-s)^2
 $$
and
 \BE
 &&\|(z_t-z_s)^*(z_t-z_s)-y(s,t)\|\\
 &&\approx \|\sum p_k^{1/2}(z_t-z_s)^*(z_t-z_s)p_k^{1/2} -y(s,t)\|\\
 && =\|\sum p_k^{1/2}\{((z_t-z_s)^*-(u_{(1-k/N)t}-u_{(1-k/N)s})^*)(z_t-z_s)\\
 && \ \ +(u_{(1-k/N)t}-u_{(1-k/N)s})^*(z_t-z_s-(u_{(1-k/N)t}-u_{(1-k/N)s}))\}p_k^{1/2}\|\\
 &&\leq 2\eps'(\|z_t-z_s\|+(C+\eps')(t-s)).
 \EE
Hence we have that
 $$
 \|z_t-z_s\|^2<(C+\eps')^2(t-s)^2+2\eps'(C+\eps')(t-s)+2\eps'\|z_t-z_s\|+\eps'',
 $$
or
 $$
 (\|z_t-z_s\|-\eps')^2<(C+\eps')^2(t-s)^2+2\eps'(C+\eps')(t-s)+\eps'{^2}+\eps''
 $$
for an arbitrarily small constant $\eps''>0$. Assuming that 
$t-s>\sqrt{2\eps'}$ and $\eps'{^2}>\eps''$, we have that
 $$
 (\|z_t-z_s\|-\eps')^2<(C+\eps'+\sqrt{\eps'})^2(t-s)^2,
 $$
which implies that
 \BE
 \|z_t-z_s\|&<&(C+\eps'+\sqrt{\eps'})(t-s)+\eps'\\
     &<&(C+\eps+2\sqrt{\eps'})(t-s)\\
     &<&(C+3\sqrt{\eps'})(t-s).
 \EE
We choose a sequence $(t_i)_{i=0}^m$ in $[0,1]$ such that 
$t_0=0<t_1<\cdots<t_m=1$ and
 $$
 \sqrt{2\eps'}<t_i-t_{i-1}<2\sqrt{\eps'}<\frac{\eps}{C+8\sqrt{\eps'}}.
 $$
Noting that $\||z_{t_i}|-1\|<3\eps'\ll\eps$, define 
$v_{t_i}=z_{t_i}|z_{t_i}|^{-1}\in \U(B_0)$. Then we have that for 
$x\in\F$, 
 \BE
 \|v_{t_i}x-xv_{t_i}\|&\leq& 
 2\|v_{t_i}-z_{t_i}\|+\|z_{t_i}x-xz_{t_i}\|\\
 &\leq& 2\|1-|z_{t_i}|\|+\|[z_{t_i},x]\|\\
 &<&6\eps'+\eps<2\eps
 \EE
and that 
 \BE
 \|v_{t_i}-v_{t_{i-1}}\|&\leq&\|v_{t_i}(1-|z_{t_i}|)\|+\|z_{t_i}-z_{t_{i-1}}\| 
 +\|v_{t_{i-1}}(1-|z_{t_{i-1}}|)\|\\
 &\leq&(C+3\sqrt{\eps'})(t_i-t_{i-1})+6\eps'.
 \EE
Thus, since $5\sqrt{2}>6$, we have that
 $$
 \|v_{t_i}-v_{t_{i-1}}\|<(C+8\sqrt{\eps'})(t_i-t_{i-1})<\eps.
 $$
Let 
 $$
 \sqrt{-1}h_i=\ln(v_{t_i}v_{t_{i-1}}^*)=\ln(1-(1-v_{t_i}v_{t_{i-1}}^*));
 $$
then
 \BE
 \|h_i\|&\leq& \|1-v_{t_i}v_{t_{i-1}}^*\|+  \|1-v_{t_i}v_{t_{i-1}}^*\|^2\\
 &<& (C+8\sqrt{\eps'})(t_i-t_{i-1})+ 
 (C+8\sqrt{\eps'})^22\sqrt{\eps'}(t_i-t_{i-1})\\
 &<& (C+(3C^2+8)\sqrt{\eps'})(t_i-t_{i-1}).
 \EE
We define, for $t\in [t_{i-1},t_i]$,
 $$
 v_t=\exp(\sqrt{-1}(\frac{t-t_{i-1}}{t_i-t_{i-1}})h_i)v_{t_{i-1}}.
 $$
Then $(v_t)$ is a path in $\U(B_0)\cap (B_0+1)$. Moreover $(v_t)$ 
is a rectifiable path satisfying that
 $$
 \length((v_y)_{y\in[s,t]})<(C+(3C^2+8)\sqrt{\eps'})(t-s).
 $$
If $t\in [t_{i-1},t_i]$, then $\|v_t-v_{t_{i-1}}\|<\eps$ and hence 
$\|\Ad\,v_t(x)-x\|<4\eps$ for $x\in\F$. Since 
$\pi(z_1)\xi_1=\eta_1$ and $\||z_t|-1\|<3\eps'$, we have that  
$\|\pi(v_1)\xi_1-\eta_1\|<3\eps'$.   This is how to construct the 
path $(v_t)$ in $\U(B_0)$ in the case 
$\lan\pi(e_{j-1})\xi,\xi\ran>\eps'{^2}$. Otherwise we simply set 
$v_t=1\in \U(B_0)$. 

Let $B_1$ be the closure of $(1-e_{j+N})B(1-e_{j+N})$. If 
$\lan(1-\pi(e_{j+N}))\xi,\xi\ran>\eps'{^2}$, set 
$f=(1-e_{j+N})^{1/2}\in B_1+1$,  
$\xi_2=\|\pi(f)\xi\|^{-1}\pi(f)\xi$, and 
$\eta_2=\|\pi(f)\eta\|^{-1}\pi(f)\eta$. There is a rectifiable 
path $(w_t)$ in $\U(B_1)$ such that $w_1=1$, 
$\pi(w_1)\xi_2=\eta_2$, and 
$\length((w_u)_{u\in[s,t]})\leq\pi(t-s)$. If 
$\lan(1-\pi(e_{j+N}))\xi,\xi\ran\leq\eps'{^2}$, set 
$w_t=1\in\U(B_1)$. Let $u_t=(v_t-1)+(w_t-1)+1\in \U(B)$. 

Note that if $\lan\pi(e_{j-1})\xi,\xi\ran>\eps'{^2}$,
 $$
 |\|\pi(e_{j-1}^{1/2})\xi\|-\|\pi(e_{j-1}^{1/2})\eta\||< 
 \delta\eps'/2<\eps'
 $$
and that if $\lan\pi(1-e_{j+N})\xi,\xi\ran>\eps'{^2}$,
 $$
 |\|\pi((1-e_{j+N})^{1/2})\xi\|-\|\pi((1-e_{j+N})^{1/2})\eta\||
 <\delta\eps'/2<\eps'.
 $$                                                                               
Note also that
 \BE
 &&\|\xi-\pi(e_{j-1}^{1/2})\xi-\pi((1-e_{j+N})^{1/2})\xi\|^2 \\
 &&= 1+\lan\pi(e_{j-1})\xi,\xi\ran+\lan(1-\pi(e_{j+N})\xi,\xi\ran
 -2\lan 
 \pi(e_{j-1}^{1/2})\xi,\xi\ran-2\lan\pi((1-e_{j+N})^{1/2})\xi,\xi\ran\\
 &&\leq \lan\pi(e_{j+N}-e_{j-1})\xi,\xi\ran<\eps'{^2}/2.
 \EE
Thus $\pi(u_1)\xi\approx\eta$ up to the order $\eps'$. Since 
$C\geq \pi$ and $\|u_t-u_s\|=\max(\|v_t-v_s\|,\|w_t-w_s\|)$, it 
follows that $\length((u_t))<C+(3C^2+8)\sqrt{\eps'}$. Since $w_t$ 
commutes with $\F$, it also follows that 
$\|\Ad\,u_t(x)-x\|<4\eps,\ x\in\F$. We then find a rectifiable 
path $(w_t)$ in $\U(B)$ such that $w_0=1$, $\pi(w_1u_1)\xi=\eta$, 
$\|w_t-1\|$ is of order $\eps'$, and $\length((w_t))$ is of order 
$\eps'$, and then form a new path connecting $(u_t)$ with 
$(w_tu_1)$, which is the desired path. This completes the proof. 
\end{pf}

\medskip
{\em Proof of Proposition \ref{R2}} \ If the hereditary 
$C^*$-subalgebra $A_1$ has a unit, this follows from (the proof 
of) \ref{A11}. If $A_1$ has an approximate identity consisting of 
projections, this also follows from \ref{A11}. Thus we assume at 
least that $A_1$ has no unit.

Let $\F$ be a finite subset of $A_1$ and $\eps>0$. We may assume 
that $\|x\|\leq1$ for $x\in\F$, $\eps>0$ is sufficiently small, 
and that there is an $e\in A_1$ such that $0\leq e\leq 1$ and 
$exe=x,\ x\in\F$. Let $B$ be a \cstar\ containing $A_1$. By the 
previous lemma we may assume that $B=\ol{A_1BA_1}$. 

Let $D$ be a \cstar\ such that $B$ is a hereditary 
$C^*$-subalgebra of $D$, $D\supset A$, and the following diagram 
is commutative:
 $$
 \begin{array}{ccc}A_1&\subset &B \\
                   \cap && \cap\\
                   A&\subset &D
  \end{array}
  $$
                                                      
Let $\eps'>0$ be such that $2\sqrt{\eps'}(C+8\sqrt{\eps'})<\eps$ 
and let $M,N\in\N$ be so large that $M\eps'{^2}>4$ and 
$N\eps'>2C$. There is a sequence $(e_0,e_1,\ldots,e_{M(N+1)})$ in 
$A_1$ such that $e_0xe_0=x,\ x\in\F$, $0\leq e_i\leq 1$, and 
$e_ie_{i-1}=e_{i-1}$. Let $d$ be a very large constant and let 
$\F_1=\F\cup \{de_i\ |\ 0\leq i\leq M(N+1)\}$. For $(\F_1,\eps)$ 
with $\F_1$ as a subset of $A$, we choose $(\G_1,\delta)$ as in 
Property \ref{R1} for $(A,C)$. We set 
 $$
 \G=\{e_{M(N+1)}xe_{M(N+1)}\ |\ x\in\G_1\}\cup\{e_i\ |\ 0\leq i\leq M(N+1)\},
 $$
which is a finite subset of $A_1$.

Let $\pi$ be an irreducible representation of $B$ and $\xi,\eta$ 
unit vectors in $\calH_{\pi}$. We extend $\pi$ to an irreducible 
representation $\rho$ of $D$; $\calH_{\pi}\subset \calH_{\rho}$. 
Suppose that
 $$
 |\lan\pi(x)\xi,\xi\ran-\lan\pi(x)\eta,\eta\ran|<\delta\eps'{^2}/2,\ \ x\in\G.
 $$
There exists a $j\in\{1,2,\ldots,(M-1)(N+1)\}$ such that
 $$
 \lan\pi(e_{j+N}-e_{j-1})\xi,\xi\ran+\lan\pi(e_{j+N}-e_{j-1})\eta,\eta\ran 
 <\eps'{^2}/2.
 $$
If $\lan\pi(e_{j-1})\xi,\xi\ran>\eps'{^2}$, then we have that for 
the normalized vectors $\xi_1,\eta_1$ for $\pi(e_{j-1}^{1/2})\xi$, 
$\pi(e_{j-1}^{1/2})\eta$ respectively, 
 $$
 |\lan\pi(x)\xi_1,\xi_1\ran   -\lan\pi(x)\eta_1,\eta_1\ran|<\delta,\ \ x\in \G,
 $$
which implies that
 $$ 
  |\lan\rho(x)\xi_1,\xi_1\ran-\lan\rho(x)\eta_1,\eta_1\ran|<\delta,\ \ x\in 
  \G_1.
 $$    
Then there is a rectifiable path $(u_t)$ in $\U(D)$ such that 
$u_0=1$, $\rho(u_1)\xi_1=\eta_1$, $\|\Ad\,u_t(x)-x\|<\eps,\ 
x\in\F_1$, and $\length ((u_y)_{y\in[s,t]})<(C+\eps')(t-s)$ for 
$0\leq s<t\leq 1$. Let $B_0$ be the closure of 
$e_{j+N-1}De_{j+N-1}$, which is a hereditary $C^*$-subalgebra of 
$B$. Then, as in the proof of the previous lemma, we obtain a 
rectifiable path $(v_t)$ in $\U(B_0)$ such that 
$v_0=1,\pi(v_1)\xi_1=\eta_1$, $v_t-1\in B_0$, 
$\|\Ad\,v_t(x)-x\|<4\eps,\ x\in\F$, and $\length((v_y)_{y\in 
[s,t]})<(C+(3C^2+8)\sqrt{\eps'})(t-s)$ (by making $d$ sufficiently 
large).  On the other hand if  
$\lan\pi(e_{j-1})\xi,\xi\ran\leq\eps'{^2}$, we set $v_t=1\in 
\U(B_0)$. 

Let $B_1$ be the closure of $(1-e_{j+N})B(1-e_{j+N})$, which is a 
hereditary $C^*$-subalgebra of $B$.  If 
 $\lan(1-\pi(e_{j+N}))\xi, \xi\ran> \eps'{^2},
 $ 
let $\xi_2$ and $\eta_2$ be the normalized vectors for 
$\pi((1-e_{j+N})^{1/2})\xi$ and $\pi((1-e_{j+N})^{1/2})\eta$ 
respectively. There is a rectifiable path $(w_t)$ in $\U(B_1)$ 
such that $w_0=1$, $\pi(w_1)\xi_2=\eta_2$, $w_t-1\in B_1$, and 
$\length((w_y)_{y\in[s,t]})\leq\pi(t-s)$. Note that 
$(w_t-1)x=0=x(w_t-1),\ x\in\F$. If 
$\lan(1-\pi(e_{j+N}))\xi,\xi\ran\leq\eps'{^2}$, we set 
$w_t=1\in\U(B_1)$.  Let $u_t=(v_t-1)+(w_t-1)+1\in \U(B)$, from 
which one can construct the desired path as in the proof of the 
previous lemma. \hfill $\square$ 

\medskip
Hence Propositions \ref{A11} and \ref{C5} can be extended as 
follows:       

\begin{rem}\label{R4}
If $\calC_C$ denotes the class of \cstars\ satisfying Property 
\ref{R1} for a constant $C\geq\pi$, then the following statements 
hold:
 \begin{enumerate}
 \item If $A$ is a non-unital \cstar, $A\in\calC_C$ 
 if and only if $\tilde{A}\in\calC_C$, where $\tilde{A}$ is the \cstar\
 obtained by adjoining a unit to $A$.
 \item If $A_1,A_2\in\calC_C$, then $A_1\oplus A_2\in\calC_C$.
 \item If $A\in\calC_C$ and $A_1$ is a hereditary $C^*$-subalgebra of $A$,
 then $A_1\in\calC_C$.
 \item If $A\in\calC_C$ and $I$ is an ideal of $A$, then $I,A/I\in\calC_C$.
 \item If $(A_n)$ is an inductive system with $A_n\in\calC_C$, 
 then $\lim_nA_n\in\calC_C$. 
 \end{enumerate}
\end{rem} 
\begin{rem}\label{R5}
If $A$ is a (unital or non-unital) \cstar\ with Property \ref{R1} 
and $C$ is a (unital or non-unital) commutative \cstar, then 
$A\otimes C$ has Property \ref{R1} for the same constant. 
\end{rem}

\medskip
\small

\end{document}